\documentclass[journal]{IEEEtran}
\renewcommand{\baselinestretch}{1.0}
\renewcommand{\theenumii}{\theenumi.\arabic{enumii}}

\usepackage[utf8]{inputenc}
\usepackage{subfigure}
\usepackage{amssymb}
\usepackage{color,soul}
\usepackage[T1]{fontenc}
\usepackage{multicol}
\usepackage{multirow}
\usepackage{textcomp}
\usepackage{multirow}
\usepackage{bigdelim}
\usepackage{graphicx}
\usepackage{epstopdf}
\usepackage{amssymb}
\usepackage{amsmath}
\usepackage{psfrag}
\usepackage{amsthm}
\usepackage{subfigure}
\usepackage{algorithmic,algorithm}
\usepackage{mathrsfs}
\usepackage{framed}
\usepackage{color}
\usepackage{multirow,enumerate}
\usepackage{bigdelim}
\definecolor{shadecolor}{rgb}{0.92,0.92,0.92}
\usepackage{soul}
\usepackage{color}
\usepackage{cuted}
\usepackage{cite}
\theoremstyle{definition}

\newtheorem{theorem}{Theorem}

\newtheorem{example}{Example}

\newtheorem{lemma}{Lemma}

\newtheorem{remark}{Remark}

\makeatletter
\newcommand{\vast}{\bBigg@{3.2}}
\newcommand{\Vast}{\bBigg@{4.5}}

\makeatother

\def\BibTeX{{\rm B\kern-.05em{\sc i\kern-.025em b}\kern-.08em
		T\kern-.1667em\lower.7ex\hbox{E}\kern-.125emX}}

\usepackage{etoolbox}
\makeatletter 
\pretocmd\@bibitem{\color{black}\csname keycolor#1\endcsname}{}{\fail}
\newcommand\citecolor[1]{\@namedef{keycolor#1}{\color{blue}}}
\makeatother

\begin{document}

\title{Data-Driven Linear Quadratic Optimization for Controller Synthesis with Structural Constraints}

\author{
	Jun Ma,  
	Zilong Cheng,  
	Xiaocong Li, 
	Wenxin Wang,\\
	Masayoshi Tomizuka, \IEEEmembership{Life Fellow,~IEEE,}
	and Tong Heng Lee	
\thanks{Jun Ma is with the Robotics and Autonomous Systems Thrust, The Hong Kong University of Science and Technology (Guangzhou), Guangzhou, China, also with the Department of Electronic and Computer Engineering, The Hong Kong University of Science and Technology, Hong Kong SAR, China, and also with the HKUST Shenzhen-Hong Kong Collaborative Innovation Research Institute, Futian, Shenzhen, China (e-mail: jun.ma@ust.hk).}
	\thanks{Zilong Cheng, Wenxin Wang, and Tong Heng Lee are with the Department of Electrical and Computer Engineering, National University of Singapore, Singapore 117583 (e-mail: zilongcheng@u.nus.edu; wenxin.wang@u.nus.edu; eleleeth@nus.edu.sg).}
		\thanks{Xiaocong Li is with the John A. Paulson School of Engineering and Applied Sciences, Harvard University, Cambridge, MA 02138 USA (e-mail: xiaocongli@seas.harvard.edu)}
	\thanks{Masayoshi Tomizuka is with the Department of Mechanical Engineering, University of California, Berkeley, CA 94720 USA (e-mail: tomizuka@berkeley.edu).}
	\thanks{This work has been submitted to the IEEE for possible publication. Copyright may be transferred without notice, after which this version may no longer be accessible.}
}

\maketitle

\begin{abstract}
For various typical cases and situations where the formulation results in an optimal control problem,
the Linear Quadratic Regulator (LQR) approach and its variants continue to be highly attractive.
In certain scenarios, it can happen that some prescribed structural constraints on the gain matrix would arise.
Consequently then, the Algebraic Riccati Equation (ARE) is no longer applicable in a straightforward way to obtain the optimal solution.
This work presents a rather effective alternative optimization approach based on
gradient projection.
The utilized gradient is obtained through a data-driven  methodology,
and then projected onto applicable constrained hyperplanes.
Essentially, this projection gradient determines
a direction of progression and computation for the gain matrix update
with a decreasing functional cost;
and then the gain matrix is further refined in an iterative framework.
With this formulation, a data-driven optimization algorithm is summarized
for controller synthesis with structural constraints.
This data-driven approach has the key advantage that it
avoids the necessity of precise modeling
which is always required in the classical model-based counterpart;
and thus the approach can additionally accommodate
various model uncertainties. Illustrative examples are also provided
in the work
to validate the theoretical results. 
\end{abstract}


\section{Introduction}

In many control systems cases that are commonly encountered,
it is certainly evident that the methodology of
optimal control is highly applicable and effective.
As is also commonly known, a specific case of the optimal control problem methodology which
is often and regularly utilized is
the Linear Quadratic Regulator (LQR) approach. In the LQR problem being considered,
it is possible to compose an
objective function which is defined as the sum of quadratic terms
with respective weightings on state variables and control variables.
Notably, the LQR problem admits the optimal solution
with several properties in terms of optimality and robustness,
and it is straightforward to determine the optimal solution
by solving the well-known Algebraic Riccati Equation (ARE)~\cite{anderson2007optimal}.


With the pervasive power of optimal control,
it is noteworthy that
a class of controller optimization problems can
essentially
be formulated as an equivalent LQR problem,
but with additional structural constraints imposed on the gain matrix~\cite{cao2009optimal}.
Most of these structural constraints come from the zero elements in the controller gain matrix,
which frequently appear
in the context of decentralized control, where the off-block diagonal elements of the composite controller gain are restricted to
be zero in the decentralized control architecture.
These can also arise from some application-specific factors,
such as controller structural restrictions
in view of the complexity in measurement and feedback~\cite{xia2020penalty}.
As is known, the ARE always results in an optimal solution
when there are
no additional structural constraints (thus known as a full matrix) for the LQR problem.
But if certain structural constraints are imposed on the gain matrix, the ARE typically cannot be straightforwardly used;
in which case the optimal solution then cannot be easily obtained by similarly straightforward analytical means~\cite{ma2021optimal}.
Indeed some of these structural constraints typically lead to
the NP-hardness of the optimization problem~\cite{fattahi2017convexity,fattahi2018transformation,wang2019system,wang2020optimal}.

\textcolor{black}{To address the structural constraints of the controller},
some direct methods in optimization can be used~\cite{sun2019survey}.
For example, the original optimal control problem can be converted
to a nonlinear constrained optimization problem
or a nonlinear programming problem~\cite{borwein2010convex}, and then various approaches can be used, such as sequential quadratic programming, interior point method, etc. In~\cite{geromel1979algorithm,geromel1982optimal}, a gradient projection method is presented to solve the constrained LQR problems with a decentralized gain matrix, which exhibits good properties such as the monotonic decreasing of the cost, the convergence of the algorithm, etc. Inspired by these works, linear quadratic constrained optimization is generalized to solve the problems with multiple equality constraints of other types, which is further used for controller synthesis in mechatronic design problems~\cite{ma2017integrated}; and the above properties can also be generalized in such settings. In~\cite{geromel1991convex,geromel1994decentralized}, the structural constraints imposed on the controller can be addressed in an alternate parameter space. In~\cite{lee2018primal}, a model-free primal-dual Q-learning algorithm is proposed to recover the optimal policy for LQR design, and it is claimed that the main results can be extended in several directions, such as input constrained optimal control design problems.
In~\cite{chertovskih2020indirect}, an indirect method is presented for solving optimal control problems with state constraints. Note that this approach is based on the application of a refined version of the maximum principle, where the continuity of the corresponding Lagrange multiplier is ensured under certain conditions.
There have also been significant progress and success reported
in approaches
to use the gradient descent method to solve the optimization problem if the gradient can be calculated effectively~\cite{han2017adaptive}.
Furthermore, the
methodology of a
projection operator is frequently used to deal with
the constraints~\cite{liu2017convergence,zhu2020projected}.
These are potentially very attractive and possibly highly effective methodologies;
though at this stage, it is certainly unclear from the available existing reported works on
how to do the projection on a specific type of constrained hyperplane considering the controller structure.
Also, an additional matter of crucial importance is that
the determination of the gradient relies on the mathematical model.
In some scenarios when the model is inaccurate (such as when there exist significant model uncertainties)
or rather difficult to be identified,
the model-based optimization approach
has to contend with this
barrier to provide the ``true'' optimal solution.
Although some robust optimization techniques have been employed
to ensure guaranteed robustness towards model uncertainties,
it nevertheless still
sacrifices optimality to a certain extent~\cite{kim2008robust,vayanos2012constraint}.
To address such a barrier, there have been some rather promising preliminary studies which utilize the methodology of data-driven approaches in control system synthesis problems~\cite{karimi2017data}. For instance, the iterative feedback tuning (IFT)~\cite{HJALMARSSON1998101,hjalmarsson1998iterative,huusom2009data} is a well-known data-driven method for controller optimization and has a wide application in automation. Some representative reinforcement learning algorithms have also been proposed to solve optimal control problems in a model-free framework, such as Q-learning~\cite{rizvi2017output,zhang2018data,lian2021robust}, adaptive dynamic programming~\cite{jiang2017robust,liu2020distributed,yang2021hamiltonian}, etc. 
	Despite the generality of reinforcement learning frameworks, they are yet to directly handle structural constraints as required by the control system. Also, in terms of the adaptive dynamic programming algorithm, it simplifies a decision by breaking it down into a sequence of decision steps over time. In~\cite{vrabie2009adaptive}, a policy iteration approach is presented to solve the LQR problem without the knowledge of the internal dynamics (i.e. the state matrix), which is summarized as an algorithm that solves the ARE effectively without knowing a certain part of the system. In~\cite{mir2019dfig}, \textcolor{black}{an adaptive optimal control scheme with guaranteed Lyapunov stability} is presented; and with a policy iteration algorithm based on persistent excitation, the nonlinear ARE is solved assuming unknown system dynamics. Nevertheless in these cases, it is not straightforward from the viewpoint of structural constraints as required by certain control systems. At this current state of these works, it still leaves an important and open question on how to properly solve the linear quadratic optimization problem with structural constraints through a data-driven approach. 


In this work, a data-driven optimization algorithm is developed
to solve the LQR problem with structural constraints.
The gradient of the objective function with respect to the gain matrix
is obtained from experiments by injecting an impulse signal to the system.
Experimental procedures are given and discussed in terms
of the feasibility of practical implementation.
The main contributions of this paper are essentially encapsulated thus in three aspects:
(1). Through a gradient projection method, the LQR problem with different structural constraints
can be solved in an iterative framework with guaranteed convergence;
(2). The determination of the gradient of the objective function
with respect to the gain matrix is formulated and obtained through a data-driven approach,
which demonstrates superiority over the model-based approach
in situations of significant model uncertainties;
(3). The projection gradient results in terms of zero-element constraints are given with rigorous mathematical derivations.
Effectively then, this work further showcases the positive and promising trend
to use the data-driven approach as a realistic and effective alternative methodology
to achieve optimal performance for industrial control systems.

The remainder of this paper is organized as follows.
In Section II, the problem statement of the LQR problem with structural constraints is provided.
Next, in Section III, the data-driven approach utilized here is presented to derive
the gradient estimation of the objective function with respect to the gain matrix.
Subsequently, the gradient projection results with zero-element constraints are presented in Section IV.
Section V then presents the data-driven optimization algorithm with some of its important properties.
In Section VI, numerical examples are given to demonstrate the applicability
and effectiveness
of the methodology of the proposed algorithm.
Finally, pertinent conclusions are drawn in Section VII.

\section{Problem Statement}
In the usual nomenclature, consider the linear time-invariant (LTI) system
\begin{equation}~\label{eq:system0}
\dot x=A x+ B u,
\end{equation}
with $x(0)=x_0$, where $x\in \mathbb{R}^n$ is the state vector, $u\in \mathbb{R}^m$ is the control input vector, $A\in \mathbb{R}^{ n \times  n}$ is the state matrix, $B\in \mathbb{R}^{ n \times m}$ is the input matrix.
As is well-known, the objective of the LQR problem is to compute a static state feedback controller
\begin{equation}~\label{eq:static}
u(K)=Kx,
\end{equation} with $K\in \mathbb{R}^{m\times n}$, which stabilizes the closed-loop system and minimizes the objective function
\begin{equation}~\label{eq:costfunction}
J(K)=\int_0^\infty \frac{1}{2} (x^T Q x+u(K)^T R u(K)) \, dt,
\end{equation}
where $Q \succeq 0$, $R \succ 0$. $(A,B)$ is assumed to be stabilizable, $(A, \sqrt Q)$ is assumed to be detectable. For the sake of simplicity, and also aligning with the standard practice in the LQR formulation, we ignore the parameter $K$ in the expression $u(K)$ in the following text.  
In this work,
\textcolor{black}{we suppose that the gain matrix is under structural constraints,}
which is expressed by $K\in \Phi$, and $\Phi$ represents
the set consisting of all the gain matrices that satisfy the prescribed structural constraints.
In particular,
the rather commonly occurring case of
zero-element constraints are to be considered in this work.

If there is a single constraint imposed on $K$, the optimization problem can be expressed in the following form:
\begin{align}
	\mathop{\text{min}}\limits_{K} \hspace{1.5mm}& J(K) \nonumber\\ \text{subject to} \hspace{1.5mm}&C(K)=0, \label{eq:opt1}
\end{align}
where $C(K)$ is a singular linear function of $K$ with $C(0)=0$. If multiple constraints are imposed on $K$, the constraints can be expressed as
$C_i(K)=0, \forall i=1,2,\ldots, N$, where $C_i(K)$ is a singular linear function of $K$ with $C_i(0)=0$, and $N$ is the number of structural constraints.

\section{Gradient Estimation}
Here, an important step in the procedure is
to initially obtain the gradient of the objective function
with respect to the gain matrix experimentally.
To attain this,
we first propose, as a conceptual step only,
to inject an impulse signal to the system.
Note that this is only a conceptual proposition
necessary at this stage to
develop several key analytical expressions
essential in the methodology.
In the development which then follows that,
we will show how the impulse signal
can be replaced in practical implementation
by a unit step signal,
and the equivalent expressions
based on application of the step signal
which are used equivalently.

First of all, we start with the case with $m=1$ to derive the gradient estimation. As an extension, the gradient estimation results for the case with $m\geq2$ are derived subsequently. The results for these two cases are summarized by Theorem~\ref{thm:gradient} and Theorem~\ref{thm:gradient_MIMO}, respectively. Without loss of generality, the initial state value can be set to zero to simplify the following results. Remarkably, whether the initial state value is zero or not will not influence the analysis, because the initial state value will be perturbed when we inject the impulse signal to the system. In this work, it is assumed that the effect of the impulse signal injection dominates the variations of the initial state variables.

For the sake of brevity, $x_{i}$, $\forall i=1,2,\ldots, n$ denotes the $i$th state variable in the state vector, and $u_{j}$, $\forall j=1,2,\ldots,m$ denotes the $j$th control input in the control input vector. Also, we define $A_+(K)=(sI-A-BK)^{-1}$.

Before we introduce Theorem~\ref{thm:gradient}, Lemma \ref{LemmaMatrixInverse} and Lemma \ref{LemmaKleiman} are given first, which are used in the proof of Theorem~\ref{thm:gradient} in the sequel.

\begin{lemma}~\cite{ma2019parameter}
	\label{LemmaMatrixInverse}
	Suppose $A, B  \in \mathbb{R}^{n\times n}$, both $A$ and $(A+B)$ are invertible, then $(A+B)^{-1}= A^{-1}- A^{-1}BA^{-1}+\mathcal{O}(\|(BA^{-1})^2\|)$.
\end{lemma}

\begin{lemma}~\cite{kleinman1966linear}
	\label{LemmaKleiman}
	Let $f(X):\mathbb{R}^{m \times n} \rightarrow \mathbb{R}$ be a function and $M(X):\mathbb{R}^{m \times n} \rightarrow \mathbb{R}^{n \times m}$ be a matrix function in terms of $X$, such that $f(X+\varepsilon \Delta X)=f(X)+\varepsilon \text{Tr} (M(X)\Delta X)$,  $\forall \Delta X \in \mathbb{R}^{m \times n}$ when $\varepsilon \rightarrow 0$, then it gives $df(X)/dX=M(X)^T$.
\end{lemma}

\begin{theorem} \label{thm:gradient}
	For the LTI system~\eqref{eq:system0} in the single-input case with the static state feedback controller~\eqref{eq:static} and the objective function~\eqref{eq:costfunction},
	if the impulse signal is injected to the system, the gradient of the objective function with respect to the gain matrix is given by
	\begin{align}
		\frac{dJ(K)}{dK}&=   \int_{0}^{\infty}x^T Q\frac{\partial x(K)}{\partial K}\,dt \nonumber\\&\qquad+\int_{0}^{\infty} u^T R\left(x^T+K\frac{\partial x(K)}{\partial K}\right)\,dt, \label{eq:gradient1}
	\end{align}
	with 
	\begin{align}
		\frac{\partial x(K)}{\partial K} &= \Bigg[\left(\frac{\partial x_1(K)}{\partial K}\right)^T \quad \left(\frac{\partial x_2(K)}{\partial K}\right)^T \quad    \cdots  \nonumber\\&  \hspace{3.8cm} \Bigg(\frac{\partial x_n(K)}{\partial K}\Bigg)^T\Bigg]^T, \label{eq:gradient2}
	\end{align}
	\begin{equation}
	\frac{\partial x_i(K)}{\partial K}=\left(A_+(K)B\delta e_i^T A_+(K)B\right)^T, \label{eq:gradient3}
	\end{equation}
	where
	$\delta$ represents the impulse signal, $e_{i}\in\mathbb{R}^n, \forall i=1,2,\ldots, n$ is the standard basis vector (the $i$th entry is one while the others are zero).
\end{theorem}

\noindent{\textbf{Proof of Theorem~\ref{thm:gradient}.}}
It is straightforward to obtain the gradient of the objective function with respect to the gain matrix, which is given by
\begin{align}
	\frac{d J(K)}{d  K}&= \int_{0}^{\infty}x^T Q\frac{\partial x(K)}{\partial K}\,dt+\int_{0}^{\infty} u^T R\frac{\partial u}{\partial K}\,dt \nonumber\\
	&=\int_{0}^{\infty}x^T Q\frac{\partial x(K)}{\partial K}\,dt \nonumber\\ &\qquad +\int_{0}^{\infty} u^T R\left(x^T+K\frac{\partial x(K)}{\partial K}\right)\,dt.
\end{align}
Then, it can be observed that to obtain the gradient of the objective function with respect to the gain matrix, the gradient of the state vector with respect to the gain matrix must be known, while the state variables and the control input can be easily captured from the experiment.

From Lemma \ref{LemmaMatrixInverse}, we have
\begin{align}
	& A_+(K+\varepsilon  \Delta K)\nonumber\\
	=&\left(sI-A-B(K+\varepsilon  \Delta K)\right)^{-1}\nonumber\\
	=&(sI-A-BK)^{-1}+(sI-A-BK)^{-1} \nonumber\\ &\qquad B \varepsilon \Delta K(sI-A-BK)^{-1} + \mathcal{O}(\cdot) \nonumber\\
	=&A_+(K)+A_+(K) B \varepsilon \Delta K A_+(K)+ \mathcal{O}(\cdot).
\end{align}

According to~\cite{callier2012linear}, $A_{+}(K+\varepsilon \Delta K)B\delta$ represents the value of the state vector when an impulse signal is injected to the system with full state measurable. Then, the standard basis vector $e_i$ can be used to extract the $i$th state variable in the state vector. Thus, we have
\begin{equation}
x_i(K)=e_i^T A_{+}(K)B\delta,
\end{equation}
then
\begin{equation}
x_i(K+\varepsilon \Delta K)=e_i^T A_{+}(K+\varepsilon \Delta K)B\delta.
\end{equation}
Similarly, we have
\begin{align}
	&x_i(K+\varepsilon \Delta K)\nonumber\\  
	=& x_i(K) + e_i^T A_+(K)B\varepsilon\Delta K A_+(K)B\delta + \mathcal{O}(\cdot) \nonumber\\
	=&x_i(K) + \textup{Tr} \left(e_i^T A_+(K)B\varepsilon\Delta K A_+(K)B\delta\right) + \mathcal{O}(\cdot) \nonumber\\
	=& x_i(K) + \varepsilon\textup{Tr} \left(A_+(K)B\delta e_i^T A_+(K)B \Delta K\right)+ \mathcal{O}(\cdot).
\end{align}
Then, from Lemma~\ref{LemmaKleiman}, it is easy to obtain \eqref{eq:gradient2} and \eqref{eq:gradient3}. This completes the proof of Theorem~\ref{thm:gradient}. ~\hfill{\qed}

\begin{remark}
	According to the definition of the Fr\'{e}chet derivative~\cite{andrews2010ricci}, the high-order terms do not influence the value the first-order derivative. Therefore, the equality in~\eqref{eq:gradient3} is established (instead of a first-order approximation), where the high-order terms can be eliminated directly.
\end{remark}

From the viewpoint of practical implementation, injecting the impulse signal to the system is impossible. To tackle this problem, the impulse signal can be replaced by the unit step signal, then \eqref{eq:gradient3} is replaced by
\begin{equation}
\frac{\partial x_i(K)}{\partial K}=\left(A_+(K)B\frac{\partial}{\partial t}   \left(\theta e_i^T A_+(K)B\right) \right)^T,
\end{equation}
where $\theta$ represents the unit step signal. Besides, as the numeric derivatives are employed in the gradient estimation, a low-pass filter can be utilized to the suppress the effect of noise. 

To realize the data-driven implementation of Theorem~\ref{thm:gradient}, the expression~\eqref{eq:gradient3} needs to be derived from the data of experiments.
Specifically, two experiments are required to determine the gradient. Experiment 1 aims to realize $\partial(\theta e_i^T A_+(K)B)/\partial t$ (equivalent to $\partial(e_i^T A_+(K)B\theta)/\partial t$), which represents the derivative of the state variables when injecting the unit step signal through the input channel of the system. In this paper, the superscripts ``$\langle1\rangle$'', ``$\langle2\rangle$'', and ``$\langle3\rangle$'' represent the Experiment 1, 2, and 3, respectively. We have $ \theta e_i^T A_+(K)B  =x_i^{\langle1\rangle}$, $\forall i=1,2,\ldots, n$, and it is easy to see that the values of all the state variables, i.e.  $x_1^{\langle1\rangle}, x_2^{\langle1\rangle}, \ldots, x_n^{\langle1\rangle}$ are required to be measured in Experiment 1. Then, $\partial(\theta e_i^T A_+(K)B)/\partial t$ can be easily derived by taking the differentiation of  $x_1^{\langle1\rangle}, x_2^{\langle1\rangle}, \ldots, x_n^{\langle1\rangle}$, and we have $\dot x_1^{\langle1\rangle}, \dot x_2^{\langle1\rangle}, \ldots, \dot x_n^{\langle1\rangle}$, respectively.

Experiment 2 realizes the remaining part of \eqref{eq:gradient3} and we need to inject $\dot x_1^{\langle1\rangle}, \dot x_2^{\langle1\rangle}, \ldots, \dot x_n^{\langle1\rangle}$ through the input channel of the system separately and measure the state variables $x_{1,1}^{\langle2\rangle}, x_{1,2}^{\langle2\rangle}, \ldots, x_{1,n}^{\langle2\rangle}$ (when injecting $\dot x_1^{\langle1\rangle}$), $x_{2,1}^{\langle2\rangle}, x_{2,2}^{\langle2\rangle}, \ldots, x_{2,n}^{\langle2\rangle}$ (when injecting $\dot x_2^{\langle1\rangle}$), $x_{n,1}^{\langle2\rangle}, x_{n,2}^{\langle2\rangle}, \ldots, x_{n,n}^{\langle2\rangle}$ (when injecting $\dot x_n^{\langle1\rangle}$), respectively. Here, define
\begin{align}
	x_1^{\langle2\rangle} &= \begin{bmatrix}
		x_{1,1}^{\langle2\rangle} & x_{1,2}^{\langle2\rangle} & \cdots & x_{1,n}^{\langle2\rangle}
	\end{bmatrix}^T, \nonumber\\
	&\hspace{2mm}\vdots\nonumber\\
	x_n^{\langle2\rangle} &= \begin{bmatrix}
		x_{n,1}^{\langle2\rangle} & x_{n,2}^{\langle2\rangle} & \cdots & x_{n,n}^{\langle2\rangle}
	\end{bmatrix}^T,~\label{eq:gr_siso}
\end{align}
we have \begin{equation}~\label{eq:before}
\frac{\partial x_i(K)}{\partial K}=\left(x_i^{\langle2\rangle}\right)^T,\quad \forall i=1,2,\ldots, n.
\end{equation}
Finally, the gradient can be easily determined by \eqref{eq:gradient1} and \eqref{eq:gradient2}.

Next, we develop
the results for the case with $m\geq2$.

\begin{theorem}~\label{thm:gradient_MIMO}
	For the LTI system~\eqref{eq:system0} in the multi-input case with the static state feedback controller~\eqref{eq:static} and the objective function~\eqref{eq:costfunction}, if the impulse signal is injected to the system through each input channel, the gradient of the objective function with respect to the gain matrix is given by
	\begin{IEEEeqnarray}{rCl}
		\frac{dJ(K)}{dK}&=&   \int_{0}^{\infty}  \sum \limits_{i=1}^n q_i x_i\frac{\partial x_i(K)}{\partial K}\,dt \nonumber\\
		&&+\int_{0}^{\infty} \sum \limits_{j=1}^m  r_j u_j \frac{\partial u_j(K)}{\partial K}\,dt,   \label{eq:gradient11}
	\end{IEEEeqnarray}
	with
	\begin{equation}
	\frac{\partial x_i(K)}{\partial K}=\left(A_+(K)B\delta_v e_i^T A_+(K)B\right)^T, \label{eq:gradient33}
	\end{equation}
	\begin{align}
		&\hspace{0mm}\frac{\partial u_j(K)}{\partial K} = \nonumber\\  &
		\hspace{2mm}\begin{aligned}
			\begin{matrix}
				&1 & \cdots  & n \\
				1\ldelim[{8}{0.1cm}&\sum\limits_{i=1}^n k_{ji}\frac{\partial x_i(K)}{\partial k_{11}} &\cdots & \sum\limits_{i=1}^n k_{ji}\frac{\partial x_i(K)}{\partial k_{1n}}&\rdelim]{8}{0.1cm}\\
				\vdots & \vdots & \vdots & \vdots\\
				j &x_1+\sum\limits_{i=1}^n k_{ji}\frac{\partial x_i(K)}{\partial k_{j1}} & \cdots & x_n+\sum\limits_{i=1}^n k_{ji}\frac{\partial x_i(K)}{\partial k_{jn}}\\
				\vdots & \vdots & \vdots & \vdots\\
				m & \sum\limits_{i=1}^n k_{ji}\frac{\partial x_i(K)}{\partial k_{m1}} & \cdots & \sum\limits_{i=1}^n k_{ji}\frac{\partial x_i(K)}{\partial k_{mn}}
			\end{matrix}\hspace{1mm},
		\end{aligned}\label{eq:gradient55}
	\end{align}
	where $q_i, \forall i=1,2,\ldots, n$ and $r_j, \forall j=1,2,\ldots, m$ represent the entries located at the $i$th row and the $i$th column of $Q$, $j$th row and $j$th column of $R$, respectively. $\delta_v\in \mathbb{R}^m$ represents a vector of the impulse signals. $e_{i}\in\mathbb{R}^n, \forall i=1,2,\ldots,n$ is the standard basis vector, $k_{ji}$ represents the entry of the matrix $K$ located at the $j$th row and $i$th column, ${\partial x_i(K)}/{\partial k_{11}}, \ldots, {\partial x_i(K)}/{\partial k_{mn}}$ are the entries located at the 1st row and 1st column, $\ldots$, the mth row and nth column of ${\partial x_i(K)}/{\partial K}$ in \eqref{eq:gradient33}.
\end{theorem}

\noindent{\textbf{Proof of Theorem~\ref{thm:gradient_MIMO}.}}
\eqref{eq:gradient11} is straightforward through matrix decomposition and the proof of \eqref{eq:gradient33} is similar to the one in Theorem~\ref{thm:gradient}. Next, we are going to complete the proof of \eqref{eq:gradient55}.

The system input is given by $u_j=\sum \limits_{i=1}^n k_{ji}x_i$. Then, we have
\begin{align}\label{eq:gradient_u}
	\frac{\partial u_j(K)}{\partial k_{fg}} &= \sum\limits_{i=1}^n\left(\frac{\partial k_{ji}}{\partial k_{fg}} x_i+k_{ji}\frac{\partial x_i(K)}{\partial k_{fg}}  \right)\nonumber\\
	&= e_j^T\left( \dfrac{\partial K}{\partial k_{fg}}x+K\dfrac{\partial x(K)}{\partial k_{fg}}\right),
\end{align}
where $e_j\in \mathbb{R}^m, \forall j=1,2,\ldots, m$ is the standard basis vector. Because $\partial K/\partial k_{fg}$ is a matrix with one entry as one located at the $f$th row and $g$th column while the other entries are zero, \eqref{eq:gradient_u} can be further simplified as
\begin{align}\label{eq:gradient_u2}
	&\dfrac{\partial u_j(K)}{\partial k_{fg}}=
	\begin{cases}
		\begin{aligned}
			&e_g^T x+\sum \limits_{i=1}^n k_{ji}\dfrac{\partial x_i(K)}{\partial k_{fg}} \quad &\textup{when} \,\, j=f,\\[1em]
			&\sum\limits_{i=1}^n k_{ji}\dfrac{\partial x_i(K)}{\partial k_{fg}}  \quad &\textup{when}\,\, j\neq f,
		\end{aligned}
	\end{cases}
\end{align}
where $e_g\in \mathbb{R}^n$ is a standard basis vector which the $g$th entry is one and all the other entries are zero. Arrange \eqref{eq:gradient_u2} into matrix form, then we have \eqref{eq:gradient55}. This completes the proof of Theorem~\ref{thm:gradient_MIMO}. ~\hfill{\qed}

\begin{remark}
	Here, it is pertinent to mention that we use the notation that is commonly utilized in the work of the IFT approach~\cite{HJALMARSSON1998101,hjalmarsson1998iterative}; and thus the operator $p = {d}/{dt}$, when utilized, denotes straightforwardly the operation $p x(t) = \dot{x}(t)$. Further thus, also straightforwardly, for a stable rational transfer function $H(s)$, the notation $x_1(t) = H(p) x_2(t)$ denotes the output signal $x_1(t)$ when an input signal $x_2(t)$ is passed through the system with the transfer function $H(s)$, which is equivalent to the expression $X_1(s) = H(s) X_2(s)$ in the transform domain (apart, of course, from an inconsequential term arising from initial conditions in the time domain which dies away exponentially fast) as long as $H(s)$ (and thus also $H(p)$) is a rational function (in $s$ and $p$ respectively) and Hurwitz in its denominator polynomial. 
	Following then this same standardized notation used commonly in the work of the IFT approach~\cite{HJALMARSSON1998101,hjalmarsson1998iterative}, our methodology (which essentially builds upon, and further extends the IFT approach) thus likewise uses this same adopted compact notation appropriately. It compactly combines the time-domain signals and the transform-domain transfer function in the above-mentioned notationally straightforward fashion without requiring further elaborate (and unnecessarily cumbersome) explicit interpretations when deriving the gradient information, for example, in Theorem~\ref{thm:gradient} and Theorem~\ref{thm:gradient_MIMO}.
\end{remark}

\begin{remark}
		Theorem~\ref{thm:gradient} and Theorem~\ref{thm:gradient_MIMO} are valid, given a feedback gain $K$ that stabilizes the system, and the system is operated in the closed loop. In this case, the optimal solution exists.
Additionally, Theorem~\ref{thm:gradient} is a particular case of Theorem~\ref{thm:gradient_MIMO} when the weighting matrices $Q$ and $R$ are diagonal. However, Theorem~\ref{thm:gradient} can also be used when the weighting matrices $Q$ and $R$ are non-diagonal, while Theorem~\ref{thm:gradient_MIMO} loses this property.
\end{remark}

Similarly, the gradient can be determined through two experiments. By replacing the impulse signal with the unit step signal, \eqref{eq:gradient33} is replaced by \begin{equation}\label{eq:gradient66}
\frac{\partial x_i(K)}{\partial K}=\left(A_+(K)B\frac{\partial}{\partial t}\left(\theta_v e_i^T A_+(K)B\right)\right)^T,
\end{equation}
where $\theta_v\in \mathbb{R}^m$ represents a vector of the unit step signals.

Experiment 1 aims to realize $\partial(\theta_v e_i^T A_+(K)B)/\partial t$. However, in the multi-input case,  $\theta_v$ is a vector (in the single-input case, $\theta$ is a scalar), so $\partial(\theta_v e_i^T A_+(K)B)/\partial t \neq \partial(e_i^T A_+(K)B\theta_v)/\partial t$. Therefore, transformations on \eqref{eq:gradient66} must be made to ensure the feasibility of the experiment. Take the first state variable $x_1$ as an example, we have $e_1^T=\begin{bmatrix}
1 & 0 & \cdots & 0
\end{bmatrix}$, then
\begin{align}\label{eq:long}
	&\frac{\partial}{\partial t}\left(\theta_v e_1^T A_+(K)B\right)\nonumber\\
	=&\frac{\partial}{\partial t}\left( \begin{bmatrix} \mathbf{1} & \mathbf{0}
	\end{bmatrix}A_+(K)B
	\textup{diag}\{\theta, \theta, \cdots, \theta\}\right),  
\end{align}
where $\mathbf 1$ and  $\mathbf 0$ represent the vector with all entries as one and the vector with all entries as zero, respectively.
In this operation, $\partial(\theta_v e_1^T A_+(K)B)/\partial t $ can be realized by injecting the unit step signal to $m$ input channels separately, taking the differentiation and pre-multiplying a matrix. \eqref{eq:long} can be written in the following form:
\begin{align}
	\frac{\partial}{\partial t}\left(\theta_v e_1^T A_+(K)B\right)
	&=\begin{bmatrix} \mathbf{1} & \mathbf{0}
	\end{bmatrix}  \begin{bmatrix}
	\dot x_{1,1}^{\langle1\rangle}   \quad  \cdots  \quad     \dot x_{1,m}^{\langle1\rangle}
\end{bmatrix},
\end{align}
where $x_{1,1}^{\langle1\rangle}= A_+(K)B  \begin{bmatrix} \theta &0 &\cdots &0 \end{bmatrix}^T$ can be measured with the unit step signal injected to the $1$st input channel only, $x_{1,m}^{\langle1\rangle}= A_+(K)B  \begin{bmatrix} 0 &\cdots &0 & \theta  \end{bmatrix}^T $ can be measured with the unit step signal injected to the $m$th input channel only. In a similar way, $ {\partial}(\theta_v e_i^T A_+(K)B)/{\partial t}$, $\forall i=1,2,\ldots,n$ can be obtained, and then define
\begin{align}
	\dot x_1^{\langle1\rangle}&=\frac{\partial}{\partial t}\left(\theta_v e_1^T A_+(K)B\right)\nonumber=\begin{bmatrix} \mathbf{1} & \mathbf{0}
	\end{bmatrix}    \begin{bmatrix}
	\dot x_{1,1}^{\langle1\rangle}   \quad  \cdots  \quad     \dot x_{1,m}^{\langle1\rangle}
\end{bmatrix} ,\nonumber\\
&\,\,\,\vdots  \nonumber \\ 
	\dot x_n^{\langle1\rangle}&=\frac{\partial}{\partial t}\left(\theta_v e_n^T A_+(K)B\right)=\begin{bmatrix} \mathbf{0} & \mathbf{1}
	\end{bmatrix} \begin{bmatrix}
	\dot x_{n,1}^{\langle1\rangle}   \quad  \cdots  \quad     \dot x_{n,m}^{\langle1\rangle}
\end{bmatrix}. ~\label{eq:xndot}
\end{align}


Note that because $\dot x_{1,1}^{\langle1\rangle} =\dot x_{2,1}^{\langle1\rangle}=\ldots=\dot x_{n,1}^{\langle1\rangle}, \dot x_{1,2}^{\langle1\rangle} =\dot x_{2,2}^{\langle1\rangle}=\ldots=\dot x_{n,2}^{\langle1\rangle}, \ldots, \dot x_{1,m}^{\langle1\rangle} =\dot x_{2,m}^{\langle1\rangle}=\ldots=\dot x_{n,m}^{\langle1\rangle}$, Experiment 1 consists of $m$ sub-experiments only.

For the sake of simplicity, define the following vectors to be used in the sequel:
\begin{align} ~\label{eq:mimmimmim}
	x_{1,1}^{\langle1\rangle}&= \begin{bmatrix}x_{1,1,1}^{\langle1\rangle} & x_{1,1,2}^{\langle1\rangle}& \cdots & x_{1,1,n}^{\langle1\rangle}
	\end{bmatrix}^T,\nonumber\\
	&\hspace{2mm} \vdots\nonumber\\
	x_{1,m}^{\langle1\rangle}&= \begin{bmatrix}
		x_{1,m,1}^{\langle1\rangle} & x_{1,m,2}^{\langle1\rangle}& \cdots & x_{1,m,n}^{\langle1\rangle}
	\end{bmatrix}^T.
\end{align}

Experiment 2 aims to realize the remaining part of \eqref{eq:gradient66}. Because \eqref{eq:xndot} can be explicitly represented by
\begin{align}
	\dot x_1^{\langle1\rangle}&= \begin{bmatrix}
		\dot x_{1,1,1}^{\langle1\rangle}  &\cdots & \dot x_{1,m,1}^{\langle1\rangle} \\
		\vdots & \ddots & \vdots\\
		\dot x_{1,1,1}^{\langle1\rangle}  &\cdots & \dot x_{1,m,1}^{\langle1\rangle}
	\end{bmatrix},\nonumber\\
	&\hspace{2mm} \vdots&\nonumber\\
	\dot x_n^{\langle1\rangle}&=\begin{bmatrix}
		\dot x_{1,1,n}^{\langle1\rangle}  &\cdots & \dot x_{1,m,n}^{\langle1\rangle} \\
		\vdots & \ddots & \vdots\\
		\dot x_{1,1,n}^{\langle1\rangle}  &\cdots & \dot x_{1,m,n}^{\langle1\rangle}
	\end{bmatrix},
\end{align}
we have
\begin{align}~\label{eq:after}
	&\frac{\partial x_1(K)}{\partial K} \nonumber\\=& \begin{bmatrix} A_+(K)B \begin{bmatrix}
			\dot x_{1,1,1}^{\langle1\rangle}   \\ \vdots \\ \dot x_{1,1,1}^{\langle1\rangle}  \end{bmatrix}  & \cdots &    A_+(K)B \begin{bmatrix}
			\dot x_{1,m,1}^{\langle1\rangle}   \\ \vdots \\ \dot x_{1,m,1}^{\langle1\rangle}  \end{bmatrix} \end{bmatrix}^T.  \nonumber\\
\end{align}
Notably, $A_+(K)B \begin{bmatrix}
\dot x_{1,1,1}^{\langle1\rangle}  & \cdots & \dot x_{1,1,1}^{\langle1\rangle}  \end{bmatrix}^T$ can be measured with $\dot x_{1,1,1}^{\langle1\rangle}$ injected to all $m$ input channels, where the measured state variables are denoted by $x_{1,1,1}^{\langle2\rangle}, \ldots, x_{1,1,n}^{\langle2\rangle}$. Similarly, $A_+(K)B \begin{bmatrix}
\dot x_{1,m,1}^{\langle1\rangle}  & \cdots & \dot x_{1,m,1}^{\langle1\rangle}  \end{bmatrix}^T$ can be measured with $\dot x_{1,m,1}^{\langle1\rangle} $ injected to all $m$ input channels, where the measured state variables are denoted by $x_{1,m,1}^{\langle2\rangle}, \cdots, x_{1,m,n}^{\langle2\rangle}$.

\begin{table*}\normalsize~\label{tab:exp}~\caption{Input signals and measured state variables required in the experimental procedures}
	\centering \resizebox{15cm}{!}{
		\begin{tabular}{|c|c|c|c|c|} 
			\hline
			& Experiments                   & Sub-experiments      & Input signals                                                                                           & Measured state variables                                                                           \\ 
			\hline
			\multirow{4}{*}{Single-input case} & Experiment 1                  & -                    & Unit step signal                                                                                        & $x_1^{\langle1\rangle}, x_2^{\langle1\rangle}, \ldots, x_n^{\langle1\rangle}$                      \\ 
			\cline{2-5}
			& \multirow{3}{*}{Experiment 2} & Sub-experiment 2.1   & $\dot x_1^{\langle1\rangle}$                                                                            & $x_{1,1}^{\langle2\rangle}, x_{1,2}^{\langle2\rangle}, \ldots, x_{1,n}^{\langle2\rangle}$          \\ 
			\cline{3-5}
			&                               & $\vdots$             & $\vdots$                                                                                                & $\vdots$                                                                                           \\ 
			\cline{3-5}
			&                               & Sub-experiment 2.n   & $\dot x_n^{\langle1\rangle}$                                                                            & $x_{n,1}^{\langle2\rangle}, x_{n,2}^{\langle2\rangle}, \ldots, x_{n,n}^{\langle2\rangle}$          \\ 
			\hline
			\multirow{11}{*}{Multi-input case} & \multirow{3}{*}{Experiment 1} & Sub-experiment 1.1   & \begin{tabular}[c]{@{}c@{}}Unit step signal to \\ the 1st input channel \end{tabular}                   & $x_{1,1,1}^{\langle1\rangle}, x_{1,1,2}^{\langle1\rangle}, \ldots, x_{1,1,n}^{\langle1\rangle}$    \\ 
			\cline{3-5}
			&                               & $\vdots$             & $\vdots$                                                                                                & $\vdots$                                                                                           \\ 
			\cline{3-5}
			&                               & Sub-experiment 1.m   & \begin{tabular}[c]{@{}c@{}}Unit step signal to\\ the $m$th input channel \end{tabular}                  & $x_{1,m,1}^{\langle1\rangle}, x_{1,m,2}^{\langle1\rangle}, \cdots, x_{1,m,n}^{\langle1\rangle}$    \\ 
			\cline{2-5}
			& \multirow{7}{*}{Experiment 2} & Sub-experiment 2.1.1 & \begin{tabular}[c]{@{}c@{}}$\dot x_{1,1,1}^{\langle1\rangle}$ to\\ all $m$ input channels \end{tabular} & $x_{1,1,1}^{\langle2\rangle}, x_{1,1,2}^{\langle2\rangle}, \ldots, x_{1,1,n}^{\langle2\rangle} $   \\ 
			\cline{3-5}
			&                               & $\vdots$             & $\vdots$                                                                                                & $\vdots$                                                                                           \\ 
			\cline{3-5}
			&                               & Sub-experiment 2.1.m & \begin{tabular}[c]{@{}c@{}}$\dot x_{1,m,1}^{\langle1\rangle}$ to\\ all $m$ input channels \end{tabular} & $x_{1,m,1}^{\langle2\rangle}, x_{1,m,2}^{\langle2\rangle}, \ldots, x_{1,m,n}^{\langle2\rangle} $   \\ 
			\cline{3-5}
			&                               & $\vdots$             & $\vdots$                                                                                                & $\vdots$                                                                                           \\ 
			\cline{3-5}
			&                               & Sub-experiment 2.n.1 & \begin{tabular}[c]{@{}c@{}}$\dot x_{1,1,n}^{\langle1\rangle}$ to\\ all $m$ input channels \end{tabular} & $x_{n,1,1}^{\langle2\rangle}, x_{n,1,2}^{\langle2\rangle}, \ldots, x_{n,1,n}^{\langle2\rangle} $   \\ 
			\cline{3-5}
			&                               & $\vdots$             & $\vdots$                                                                                                & $\vdots$                                                                                           \\ 
			\cline{3-5}
			&                               & Sub-experiment 2.n.m & \begin{tabular}[c]{@{}c@{}}$\dot x_{1,m,n}^{\langle1\rangle}$ to\\ all $m$ input channels \end{tabular} & $x_{n,m,1}^{\langle2\rangle}, x_{n,m,2}^{\langle2\rangle}, \ldots, x_{n,m,n}^{\langle2\rangle} $   \\ 
			\cline{2-5}
			& Experiment 3                  & -                    & \begin{tabular}[c]{@{}c@{}}Unit step signal to \\ all $m$ input channels \end{tabular}                  & $x_1^{\langle3\rangle}, x_2^{\langle3\rangle}, \ldots, x_n^{\langle3\rangle}$                      \\
			\hline
		\end{tabular}
}
\end{table*}

Define the following vectors:
\begin{align}~\label{eq:mim}
	x_{1,1}^{\langle2\rangle} &=  \begin{bmatrix}
		x_{1,1,1}^{\langle2\rangle}  & x_{1,1,2}^{\langle2\rangle} &\cdots & x_{1,1,n}^{\langle2\rangle} \end{bmatrix}^T, \nonumber\\
	&\hspace{2mm} \vdots\nonumber\\
	x_{1,m}^{\langle2\rangle} &=  \begin{bmatrix}
		x_{1,m,1}^{\langle2\rangle}  & x_{1,m,2}^{\langle2\rangle} &\cdots & x_{1,m,n}^{\langle2\rangle} \end{bmatrix}^T,
\end{align}
then, \eqref{eq:after} is expressed by
\begin{equation}~\label{eq:gr}
\frac{\partial x_1(K)}{\partial K} = \begin{bmatrix}   x_{1,1}^{\langle2\rangle} & x_{1,2}^{\langle2\rangle} & \cdots & x_{1,m}^{\langle2\rangle} \end{bmatrix}^T.
\end{equation}
Similarly, all the other state variables can be measured with the same technique as mentioned above.

Experiment 3 is required to measure $x_1, x_2, \ldots, x_n$, this step can be easily done by injecting the unit step signal to all $m$ input channels and take the differentiation of the state variables, i.e. $x_i^{\langle3\rangle}=\partial({e_i^T A_+(K)B \theta_v)}/\partial t$, $\forall i=1,2,\ldots,n$. Then, the gradient can be easily determined by \eqref{eq:gradient11} and \eqref{eq:gradient55}.

To summarize the input signals and measured state variables in each experiment and sub-experiment, Table 1 is shown.

	\begin{remark}
		Theorem~\ref{thm:gradient} and Theorem~\ref{thm:gradient_MIMO} demonstrate the results with the impulse signal injected to the system.
		However, as also noted and indicated earlier, it is
		the unit step signal which will be used in practical implementation
		(with the additionally developed equivalent analytical relationships).
		It is also essential
		to carefully note that in the developments,
		the measured system state variables and input signals (which are denoted by
		respective symbols with superscripts ``$\langle1\rangle$'', ``$\langle2\rangle$'', and ``$\langle3\rangle$'')
		are different from the
		ones in Theorem~\ref{thm:gradient} and Theorem~\ref{thm:gradient_MIMO}
		which are annotated differently;
		and that the methodology utilizes
		numerical derivatives.
	\end{remark}
	
	\section{Gradient Projection}
	On the basis of~\cite{geromel1982optimal}, the gradient projection results with zero-element constraints are presented in the following contents.
	Assume $dJ(K)/dK$ is not null, then we aim to
	determine the projection gradient matrix $D$ which guarantees the  decrease of the objective function with a step size $\alpha$, i.e. $ J(K-\alpha D) < J(K)$ with $\alpha >0$.
	The Frobenius norm of a matrix is denoted by $\|\cdot\|_F$. The optimization problem~\eqref{eq:opt1} is equivalent to calculate $D$ such that the Euclidean distance between $D$ and $dJ(K)/dK$ is minimized.
	Therefore, for the problem with a single constraint, it is straightforward that the optimization problem \eqref{eq:opt1} can be solved by a gradient descent approach, where the constrained gradient matrix can be determined by
	\begin{align} 
		\mathop{\text{min}}\limits_{D} \hspace{1.5mm}&\frac{1}{2} \left\|\frac {dJ(K)}{d{K}}-D\right\|_F^2 \nonumber\\
		\text{subject to} \hspace{1.5mm}&C(D)=0. \label{eq:opt2}
	\end{align}
	The dual problem of \eqref{eq:opt2} is given by
	\begin{equation}
	\mathop{\text{max}} \limits_{\Lambda} \mathop{\text{min}}\limits_{D} \hspace{1.5mm} \left( \frac{1}{2} \left\|\frac {dJ(K)}{d{K}}-D \right\|_F^2 + \text{Tr}\left(\Lambda^T C(D)\right) \right),   \label{eq:opt3}
	\end{equation}
	where $\Lambda$ is the Lagrange multiplier with respect to the structural constraint $C(D)=0$. For the problem with multiple constraints, its dual problem can be expressed similarly, 
	where $\Lambda_i$ is the Lagrange multiplier associated with the structural constraint $C_i(D)=0$. Remarkably, we assume that the primal problem \eqref{eq:opt2} is strictly feasible. Then, Slater’s condition is satisfied and strong duality always holds. Thus, the optimal solution can be determined by solving the dual problem \eqref{eq:opt3}.
	In this work, we consider the structural constraints where some elements in the controller gain matrix are zero, which can be represented by $C(D)=GDH$ (for a single constraint) or $C_i(D)=G_iDH_i, \forall i=1,2,\ldots,N$ (for multiple constraints), respectively. Remarkably, $G$, $G_i$, $G_j$ are full row rank matrices and $H$, $H_i$, $H_j$ are full column rank matrices. These matrices are required to be chosen depending on the specific structure of the controller gain matrix.

	Note that in~\cite{geromel1982optimal}, the determination of projected gradient is given considering a single linear constraint, and	Lemma~\ref{thm:projection_zero_single} 
	restates part of the results in~\cite{geromel1982optimal} and generalizes the determination of the projection gradient in two scenarios considering the projection onto a single constrained hyperplane and multiple constrained hyperplanes.
	
	\begin{lemma}\label{thm:projection_zero_single}
		The following statements hold:\\
	(i) The gradient projected onto a single constrained hyperplane $C(D)= GDH$ is given by
			 \begin{IEEEeqnarray*}{rCl}
				D=\frac {dJ(K)}{d{K}}
				- G^T \left(GG^T\right)^{-1}G \frac{dJ(K)}{d{K}} H \left(H^T H\right)^{-1}  H^T; \\ \yesnumber \label{eq:projection_gradient_zero_single}
			\end{IEEEeqnarray*}
	(ii) The gradient projected onto multiple constrained hyperplanes $C_i(D)= G_iDH_i$, $\forall i=1,2,\ldots, N$ is given by
			\begin{equation}
			D=\frac {dJ(K)}{d{K}}- \sum_{i=1}^N G_{i}^T \Lambda_i  H_{i}^T,
			\end{equation}
			where $\Lambda_i, \forall i=1,2,\ldots,N$ can be obtained by solving the following equations:
			\begin{align}
				G_1 \frac {dJ(K)}{d{K}} H_1 &=   \sum_{i=1}^N G_1 G_i^T \Lambda_i H_i^T H_1 , \nonumber\\
				&\,\,\,\vdots \nonumber\\
				G_N \frac {dJ(K)}{d{K}} H_N &=   \sum_{i=1}^N G_N G_i^T \Lambda_i H_i^T H_N.
			\end{align}

	\end{lemma}
	
	\noindent{\textbf{Proof of Lemma~\ref{thm:projection_zero_single}.}}
	Part of the proof has been provided in~\cite{geromel1982optimal}. For Statement $(i)$, the necessary and sufficient optimal solution to the dual problem \eqref{eq:opt3} is given by
	\begin{equation}
	D-\frac {dJ(K)}{d{K}}+\frac {\partial}{\partial {D}}\left(\text{Tr}\left(\Lambda^T C(D)\right)\right) =0,
	\end{equation}
	which gives
	\begin{equation}
	D=\frac {dJ(K)}{d{K}}- G^T  \Lambda  H^T. \label{eq:projection1_zero_single}
	\end{equation}
	Since $GDH=0$, it is easy to see that
	\begin{equation}
	G\frac {dJ(K)}{d{K}}H- GG^T  \Lambda  H^TH=0.
	\end{equation}
	Therefore, we have
	\begin{equation}
	\Lambda = \left(GG^T\right)^{-1}G \frac{dJ(K)}{d{K}} H \left(H^T H\right)^{-1}. \label{eq:projection2_zero_single}
	\end{equation}
	Substitute \eqref{eq:projection2_zero_single} to \eqref{eq:projection1_zero_single}, we have \eqref{eq:projection_gradient_zero_single}. This concludes the proof of Statement $(i)$. The proof of Statement $(ii)$ is omitted because it is similar to the proof of Statement $(i)$.~\hfill{\qed}
	
	\section{Proposed Algorithm}

\begin{figure}[h]
	\centering
	\includegraphics[trim=0 0 0 0,width=0.9\columnwidth]{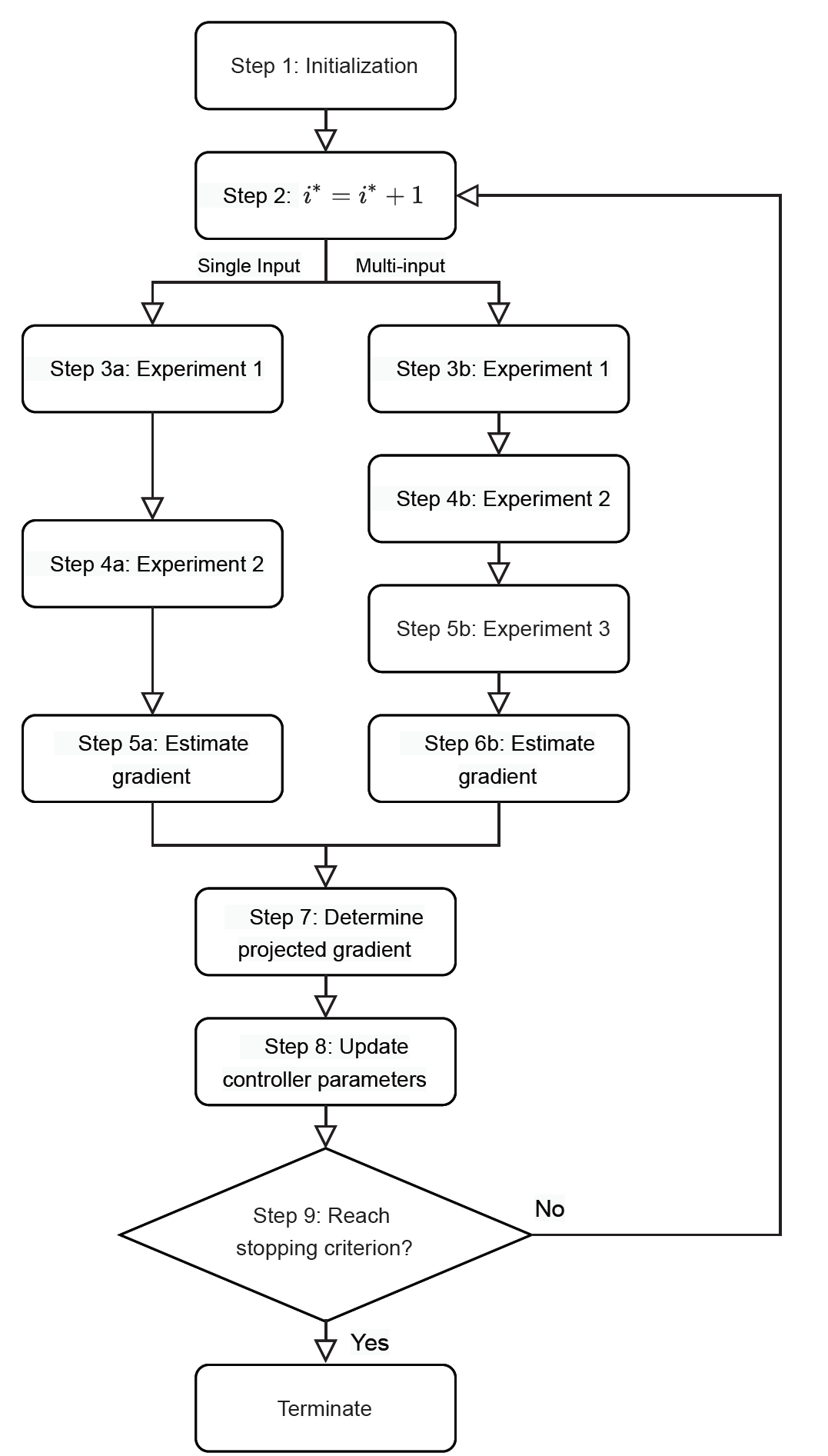}
	\caption{Flowchart of the proposed data-driven linear quadratic optimization algorithm for controller synthesis with structural constraints.}
	\label{fig:flow}
\end{figure}

	To summarize the above results on solving the controller synthesis problem under structural constraints,
	the proposed algorithm is shown in the listed Algorithm 1. Also, to ensure the clarity and conciseness of the algorithm, a flowchart is given in~ Fig.~\ref{fig:flow}.
	
	\begin{algorithm*}[!t]\label{Algo:Databased}
		\caption{Data-driven Linear Quadratic Optimization Algorithm for Controller Synthesis with Structural Constraints}
		\label{algorithm}
		\begin{itemize}
			\item \textbf{{Step 1:}}   Set the iteration number $i^*=0$, and initialize the controller $K_0$ such that the closed-loop system is stable.
			\item \textbf{{Step 2:}} $i^*=i^*+1$. For the single-input case, go to Step 3a; for the multi-input case, go to Step 3b.

			\item \textbf{{Step 3a:}} Execute Experiment 1: inject the unit step signal to the system, measure the state variables $x_1^{\langle1\rangle}, x_2^{\langle1\rangle},$ $\ldots, x_n^{\langle1\rangle}$ and the system input $u^{\langle1\rangle}$. Calculate their derivatives $\dot x_1^{\langle1\rangle}, \dot x_2^{\langle1\rangle}, \ldots, \dot x_n^{\langle1\rangle}$ and $\dot u^{\langle1\rangle}$.
			
			\item \textbf{{Step 4a:}} Execute Experiment 2: inject $\dot x_1^{\langle1\rangle}$ to the system, measure the state variables $x_{1,1}^{\langle2\rangle}, x_{1,2}^{\langle2\rangle}, \ldots, x_{1,n}^{\langle2\rangle}$; inject $\dot x_2^{\langle1\rangle}$ to the system, measure the state variables $x_{2,1}^{\langle2\rangle}, x_{2,2}^{\langle2\rangle}, \ldots, x_{2,n}^{\langle2\rangle}$; $\ldots$ ; inject $\dot x_n^{\langle1\rangle}$ to the system, measure the state variables $x_{n,1}^{\langle2\rangle}, x_{n,2}^{\langle2\rangle}, \ldots, x_{n,n}^{\langle2\rangle}$.

			\item \textbf{{Step 5a:}} Calculate  $ \partial x_i(K)/\partial K$ by~\eqref{eq:gr_siso} and \eqref{eq:before} and the estimation of ${dJ(K_{i^*})}/{dK}$ by \eqref{eq:gradient1} and \eqref{eq:gradient2}. Go to Step 7.

			\item \textbf{{Step 3b:} }Execute Experiment 1:  inject the unit step signal to the 1st input channel, measure the state variables  $x_{1,1,1}^{\langle1\rangle},$ $x_{1,1,2}^{\langle1\rangle},$ $\ldots,$ $x_{1,1,n}^{\langle1\rangle}$; inject the unit step signal to the 2nd input channel, measure the state variables  $x_{1,2,1}^{\langle1\rangle},$ $x_{1,2,2}^{\langle1\rangle},$ $\ldots,$ $x_{1,2,n}^{\langle1\rangle}$; $\ldots$ ;  inject the unit step signal to the $m$th input channel, measure the state variables  $x_{1,m,1}^{\langle1\rangle}, x_{1,m,2}^{\langle1\rangle}, \ldots, x_{1,m,n}^{\langle1\rangle}$. Calculate their derivatives $\dot x_{1,1,1}^{\langle1\rangle},$ $\dot x_{1,1,2}^{\langle1\rangle},$ $\ldots,$ $\dot x_{1,1,n}^{\langle1\rangle},$ $\dot x_{1,2,1}^{\langle1\rangle},$ $\dot x_{1,2,2}^{\langle1\rangle},$ $\ldots,$ $\dot x_{1,2,n}^{\langle1\rangle},$ $\ldots,$ $\dot x_{1,m,1}^{\langle1\rangle},$ $\dot x_{1,m,2}^{\langle1\rangle},$ $\ldots,$ $\dot x_{1,m,n}^{\langle1\rangle}$.

			\item \textbf{{Step 4b:}} Execute Experiment 2: inject $\dot x_{1,1,1}^{\langle1\rangle}$ to all $m$ input channels, measure the state variables $x_{1,1,1}^{\langle2\rangle}, x_{1,1,2}^{\langle2\rangle}, \ldots,$ $x_{1,1,n}^{\langle2\rangle} $; inject $\dot x_{1,2,1}^{\langle1\rangle}$ to all $m$ input channels, measure  the state variables $x_{1,2,1}^{\langle2\rangle}, x_{1,2,2}^{\langle2\rangle}, \ldots, x_{1,2,n}^{\langle2\rangle} $; $\ldots$ ; inject $\dot x_{1,m,1}^{\langle1\rangle}$ to all $m$ input channels, measure the state variables $x_{1,m,1}^{\langle2\rangle}, x_{1,m,2}^{\langle2\rangle}, \ldots, x_{1,m,n}^{\langle2\rangle} $; $\ldots$ ; inject $\dot x_{1,1,n}^{\langle1\rangle}$ to all $m$ input channels, measure the state variables $x_{n,1,1}^{\langle2\rangle}, x_{n,1,2}^{\langle2\rangle}, \ldots, x_{n,1,n}^{\langle2\rangle} $; inject $\dot x_{1,2,n}^{\langle1\rangle}$ to all $m$ input channels, measure  the state variables $x_{n,2,1}^{\langle2\rangle}, x_{n,2,2}^{\langle2\rangle}, \ldots, x_{n,2,n}^{\langle2\rangle} $; $\ldots$ ; inject $\dot x_{1,m,n}^{\langle1\rangle}$ to all $m$ input channels, measure the state variables $x_{n,m,1}^{\langle2\rangle}, x_{n,m,2}^{\langle2\rangle}, \ldots, x_{n,m,n}^{\langle2\rangle}$.

			\item \textbf{{Step 5b:}} Execute Experiment 3: inject the unit step signal to all $m$ input channels, measure the state variables $x_1^{\langle3\rangle},$ $x_2^{\langle3\rangle}, \ldots, x_n^{\langle3\rangle}$ and system inputs $u^{\langle3\rangle}$. Calculate their derivatives $\dot x_1^{\langle3\rangle}, \dot x_2^{\langle3\rangle}, \ldots, \dot x_n^{\langle3\rangle}$ and $\dot u^{\langle3\rangle}$.
			
			\item \textbf{{Step 6b:}} Calculate  $ \partial x_1(K)/\partial K$ by \eqref{eq:mim} and \eqref{eq:gr}, and similarly calculate $ \partial x_2(K)/\partial K, \ldots, \partial x_n(K)/\partial K$, determine the estimation of ${dJ(K_{i^*})}/{dK}$ by \eqref{eq:gradient11} and \eqref{eq:gradient55}. Go to Step 7.

			\item \textbf{{Step 7:}} Determine the projected gradient $D_{i^*}$ by Lemma~\ref{thm:projection_zero_single}.
			
			\item \textbf{{Step 8:}} Update the controller parameters by ${K}_{i^*}={K}_{(i-1)^*}-\alpha D_{i^*}$.
			
			\item \textbf{{Step 9:}} If the stopping criterion is reached, terminate the optimization process; otherwise, go back to Step 2.
			
		\end{itemize}
	\end{algorithm*}

It is worthwhile to mention that a feasible step size $\alpha$ should be specified by the user appropriately. Generally, for model-based optimization problems with the use of gradient descent technique, the step size could be selected by line search. Similarly, for data-driven optimization problems, line search can also be effectively used.
 The stopping criterion of the optimization algorithm can be implemented in different ways. One way is to define a small  $\epsilon$, such that the optimization algorithm is terminated as soon as $\|D_{i^*}\|_F\leq \epsilon$, which will be discussed later; alternatives of termination condition can be defined in terms of iteration numbers, etc.
	
	It is worthwhile to mention that the algorithm exhibits some important properties, as summarized in Theorem~\ref{property}, which is restated and generalized from \cite{geromel1982optimal}.  For clarification of Theorem~\ref{property}, the iteration number is denoted by $N^*$. Note that the evaluation of the cost till $t=\infty$ is not possible, and thus we propose the approximation of this value in a more practical way. Since the cost function value will not be influenced much after the state variables and control input variables converge to zero, it is practical to estimate the cost function value depending on the time instance when state variables and control input variables converge. 
	
	\begin{theorem}\label{property}
		Given  $K_{0} \in \mathbb{R}^{m\times n}$ such that $C({K}_{0})=0$ and $\textup{max}(\operatorname{Re}(\textup{eig}(A+B{K}_{0})))<0$, the following statements hold:\\
		(i)	The optimization algorithm ensures local convergence with $\epsilon$-optimality, i.e. $\|D_{N^*}\|_F\leq \epsilon$;\\
		 (ii) $C(K_{i^*})=0$, $\forall i^*=1,2,\ldots,N^*$;\\
		(iii) $\textup{max}(\operatorname{Re}(\textup{eig}(A+BK_{i^*})))<0$, $\forall i^*=1,2,\ldots,N^*$.
	\end{theorem}
	
	\noindent \textbf{Proof of Theorem~\ref{property}:}
	Part of the proof has been shown in~\cite{geromel1982optimal}. For Statement $(i)$,  Lemma \ref{LemmaKleiman} is used again, then we have
	\begin{equation}
	J(K-\alpha D)=J(K)-\alpha \text{Tr} \left(\left(\frac {dJ(K)}{dK}\right)^T D\right)+ \mathcal{O}(\cdot). \label{eq:pro1}
	\end{equation}
	
	Since $D=0$ is always feasible in~\eqref{eq:opt2},
	we have
	\begin{equation}
	\left\|\frac {dJ(K)}{dK}-D\right\|_F^2\leq \left\|\frac {dJ(K)}{dK}\right\|_F^2,
	\end{equation}
	which leads to
	\begin{equation}
	\text{Tr}\left(\left(\frac {dJ(K)}{dK}\right)^T D\right) \geq \frac{1}{2} \left\|D\right\|_F^2.
	\end{equation}
	From \eqref{eq:pro1}, it is easy to derive
	\begin{equation}
	J(K-\alpha D)-J(K)\leq -\frac{1}{2}\alpha \left\|D\right\|_F^2,
	\end{equation}
	with $\|D\|_F\neq 0$. Hence, if $\|D\|_F\neq 0, \exists  \bar\alpha>0$, such that $J(K-\alpha D)<J(K), \forall 0< \alpha \leq  \bar\alpha$.
	
	For the case with a single constraint, in view of the necessary optimal condition for \eqref{eq:opt2}, we have
	\begin{equation}
	\frac {dJ(K)}{dK} +G^T  \Gamma H^T =0, \label{eq:KKT1}
	\end{equation}
	with $K \in \Phi$, i.e. $G K H=0$.  Similarly, for the necessary optimal condition for \eqref{eq:opt3}, we have
	\begin{align}
		D-\frac {dJ(K)}{dK} +G^T  \Lambda H^T=0, \label{eq:KKT2}
	\end{align}
	with $D \in \Phi$, i.e. $G D H=0$.  From \eqref{eq:KKT1} and \eqref{eq:KKT2}, it is easy to observe that if $\|D\|_F=0$, the original problem is equivalent to the projection problem with $\Gamma=-\Lambda$. Hence, if $\|D\|_F= 0$, KKT conditions of optimization problem \eqref{eq:opt2} are equivalent to those of \eqref{eq:opt3}. Similar results apply to the case with multiple constraints, 
	where the algorithm locally converges to $\epsilon$-optimality in finite time.
	For Statement $(ii)$, the projection gradient satisfies the structural constraint in each iteration, i.e. $C(D_{i^*})=0$. Hence, $C(\alpha D_{i^*})=0$. Since ${K}_{i^*}={K}_{i^*-1}-\alpha D_{i^*}$, if $C({K}_{i^*-1})=0$, then $C({K}_{i^*})=0$. Therefore, if $\exists K_{0} \in \mathbb{R}^{m\times n}$ such that $C({K}_{0})=0$, then $C({K}_{i^*})=0$, $\forall i^*=1,2, \ldots ,N^*$.
	For Statement $(iii)$, since any stable closed-loop system corresponds to a finite cost and the cost is monotonically decreasing during iterations according to Statement $(i)$, it is trivial to conclude that the closed-loop stability is preserved as long as any initial control gain stabilizes the closed-loop system. This completes the proof of Theorem~\ref{property}. \hfill{\qed}
	
		\begin{remark}
Theoretically, with more data, the performance of the proposed algorithm is better because the estimation of the gradient becomes more accurate. However, more data generally leads to higher computational effort. In this sense, the user needs to balance the performance and computational effort by choosing the amount of data appropriately.
\end{remark}
	\begin{remark}
			In the proposed algorithm, the gradient of the state vector with respect to the gain matrix can be measured directly. Thus, we need to consider the complexity of calculating the gradient of the objective function with respect to the gain matrix and also the complexity of gradient projection. In the single input case, the gradient of the objective function with respect to the gain matrix is obtained from \eqref{eq:gradient1} and \eqref{eq:gradient2}, and the complexity is $n^4$ (when $Q$ is diagonal, the complexity is reduced to $n^2$). In the multi-input case, the gradient of the objective function  with respect to the gain matrix is obtained from \eqref{eq:gradient11} and \eqref{eq:gradient55}. Since it is assumed that $Q$ is diagonal in the multi-input case, it shows that the complexity is $n^2m+m^2n=mn(m+n)$. As for the gradient projection, it is straightforward to see that the complexity is $n^2m+m^2n = mn(m+n)$ to project onto a single constrained hyperplane, and the complexity is $Nmn(m+n)$ to project onto $N$ constrained hyperplanes.	
	\end{remark}

	\section{Illustrative Example}
	To illustrate the applicability and effectiveness of the above results,
	two examples are reproduced from \cite{geromel1982optimal} and \cite{davison1990benchmark}
	with modifications.
	First of all, Example~\ref{exam:1} presents a controller design problem
	with a single input and a single zero-element constraint.
	Next, Example~\ref{exam:2} is investigated, which presents a controller design problem
	with multiple inputs and multiple zero-element constraints.	Essentially for the control problems in \cite{geromel1982optimal} and \cite{davison1990benchmark}, the system model is assumed to be known, and this is the major difference between our approach and model-based counterpart. Although the models in the two examples are given below, these models are only used for generating the data to execute the proposed data-driven algorithm.
It is worth describing further here
	that in our detailed experimentations,
	we have indeed proceeded with explorations involving
	experimental computations with
	the proposed algorithm 
	to control a suitably larger dimension quadrotor-type system (which has 12 state variables and 4 control inputs).
	With the huge amount of data when
	the dimension of the system (particularly the dimension of control input) increases,
	and for this aspect which is essentially matters pertaining to practical realization and implementable computations,
	we can state that with the typical usual computational power which is not of the so-termed ``super-computing'' capability,
	this huge amount of data practically leads to serious computational burden 
	and, at this stage, cannot yield a practical solution. 
	Nevertheless here, it should be noted that this is an aspect which relates 
	just to matters pertaining to practical realization and implementable computations.
	In so far as
	the key theoretical properties of the methodology are concerned,
	these had been
	rigorously established in Theorems 1, 2, and 3.
	Appropriately here thus, 
	this aligns with Remarks 5 and 6 in the manuscript. 
	Considering this limitation, 
	and with current practical computational resources,
	indeed
	we can only deal with systems with smaller dimensions at this stage.

	\begin{example}\label{exam:1}
		Consider a linear 2nd-order system, where	\begin{align}\nonumber
			\dot x&=Ax+Bu, \quad x(0)=x_0,\\   u&=Kx ,\nonumber
		\end{align} and
		\begin{equation}\nonumber
		\  A=\begin{bmatrix} 0 & 1\\-1 & -\sqrt{2} \end{bmatrix}, \quad B=	\begin{bmatrix} 0 \\ 1 \end{bmatrix}. \end{equation}
		Determine an optimal controller   \begin{equation}\nonumber K=\begin{bmatrix} k_1 & k_2 \end{bmatrix},\end{equation}  with $k_1=0$, such that
		\begin{equation}\nonumber
		\mathop{\textup{min}}\limits_{K} \hspace{1.5mm} \int_0^\infty \left(x^TQx+u^TRu \right) \,dt,
		\end{equation}
		where  $Q=\begin{bmatrix} 1 & 0\\0 & 1 \end{bmatrix}$, $R=0.1$.
	\end{example}

\begin{figure}
		\centering
		\includegraphics[trim=0 0 0 0,width=0.9\columnwidth]{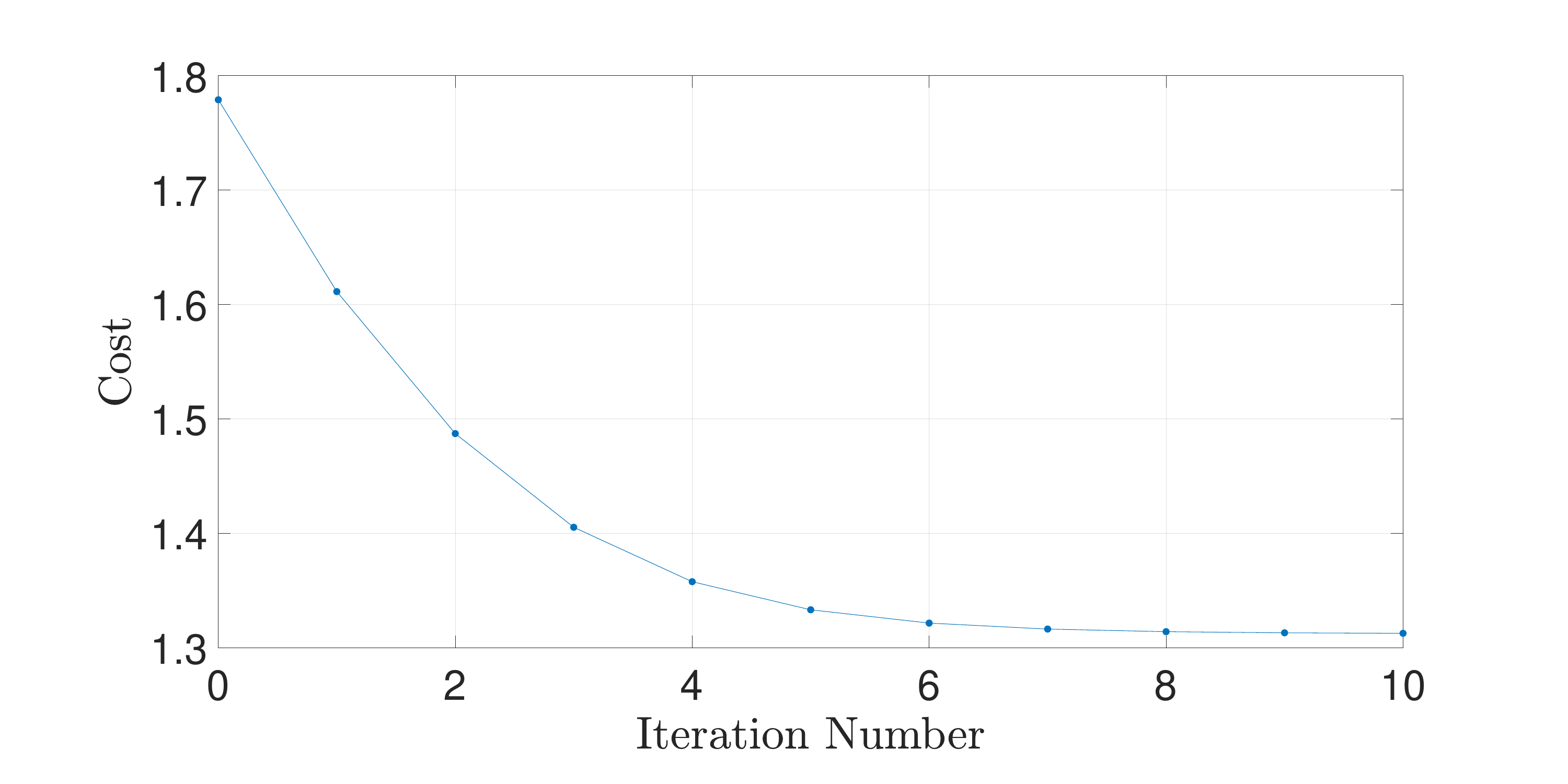}
		\caption{Change of the cost during iterations in Example 1.}
		\label{fig:cost1}
	\end{figure}

\begin{figure}
		\centering
		\includegraphics[trim=0 0 0 0,width=0.9\columnwidth]{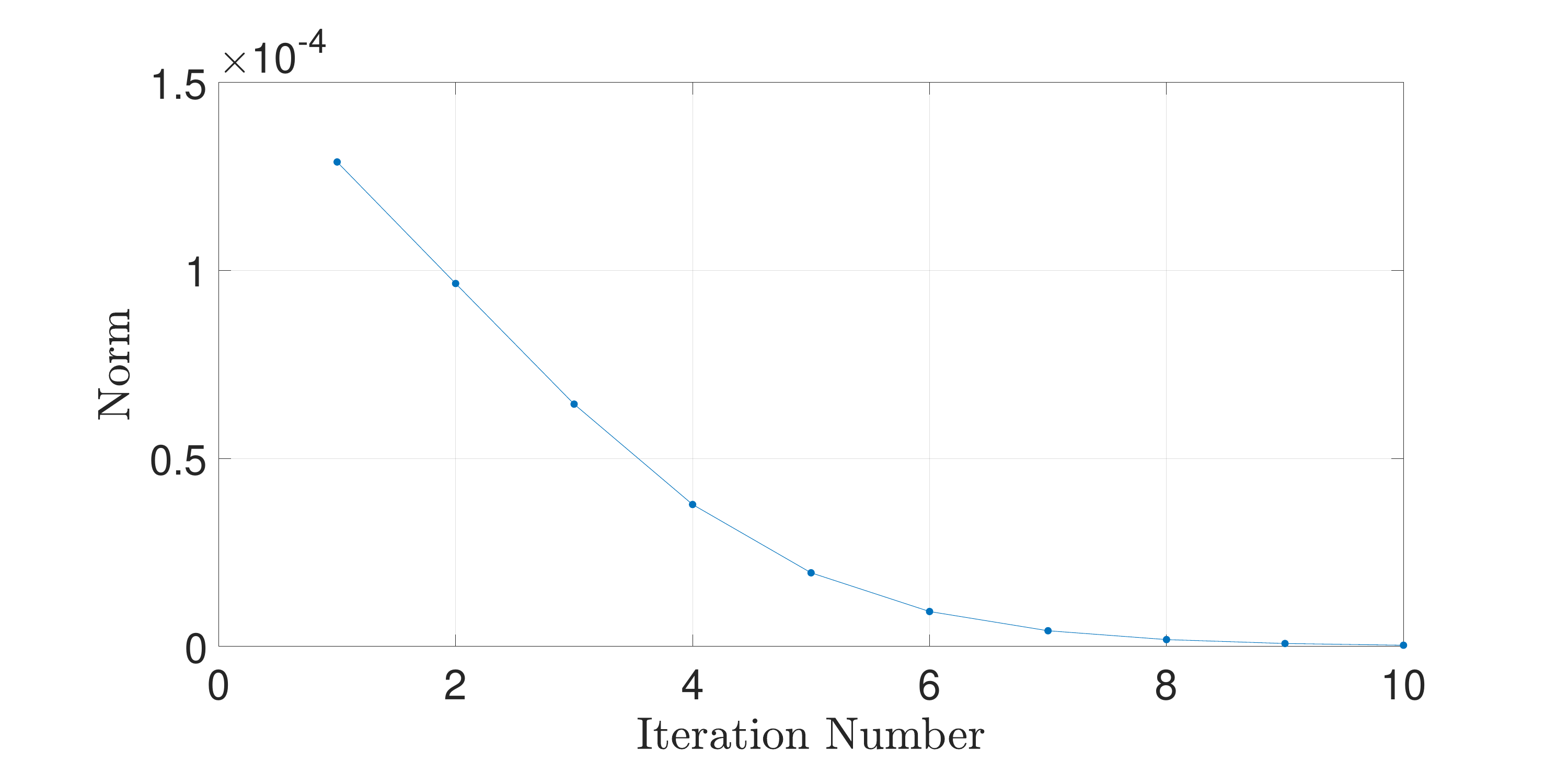}
		\caption{Change of the norm during iterations in Example 1.}
		\label{fig:norm1}
	\end{figure}

\begin{figure}
		\centering
		\includegraphics[trim=0 0 0 0,width=0.9\columnwidth]{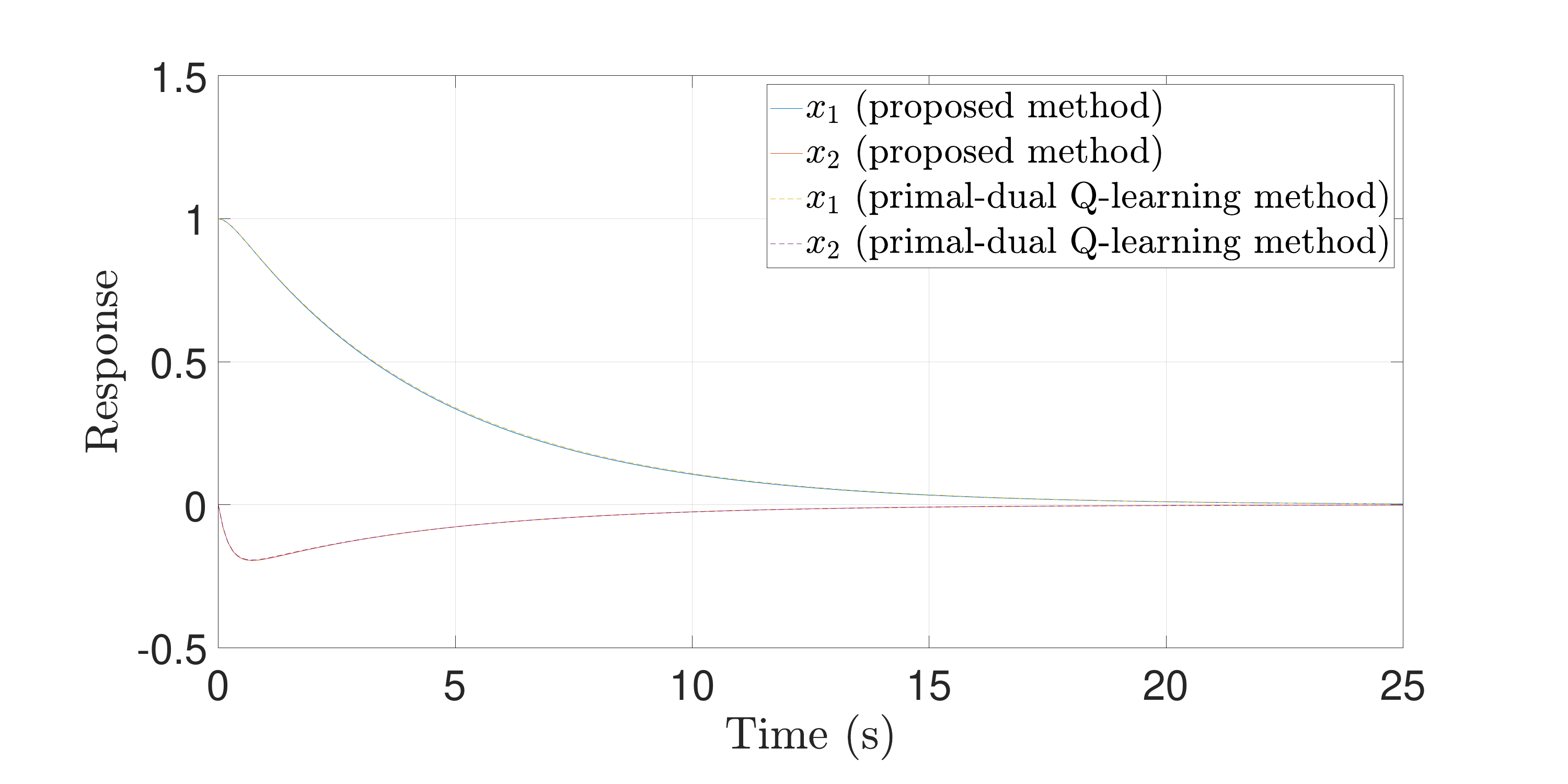}
		\caption{\textcolor{black}{System response of the proposed method and the comparative method in Example 1.}}
		\label{fig:response1}
	\end{figure}

	The structural constraint on $K$ can be stated as
	\begin{equation}\nonumber KH=0, \quad H=\begin{bmatrix} 1 & 0 \end{bmatrix}^T.\end{equation} 
	The controller gain is initialized as  \begin{equation} \nonumber K_0=\begin{bmatrix} 0&-5 \end{bmatrix},\end{equation} which stabilizes the closed-loop system with a cost of $1.7789$. Then, the optimization is executed with a stopping criteria defined as a total of 10 iterations. Fig.~\ref{fig:cost1} and Fig.~\ref{fig:norm1} illustrate the changes of the cost and the norm during 10 iterations. It can be observed that the cost and the norm drop most significantly within the first iteration. For the following iterations, the cost and the norm still decrease, but with smaller amplitudes. After 10 iterations, the cost is given by $1.3127$ with a norm of $3.429 \times 10^{-7}$, and the optimal controller gain is given by \begin{equation} \nonumber K^*=\begin{bmatrix}  0  &  -3.1813\end{bmatrix}.\end{equation} With $x_0=\begin{bmatrix}
	1&0
	\end{bmatrix}^T$ and the optimal gain matrix, the system response is shown as the solid lines in Fig.~\ref{fig:response1}.
	
	\textcolor{black}{To demonstrate the effectiveness of the proposed method, a comparison is performed using one of the representative model-free linear quadratic optimization method called the primal-dual Q-learning algorithm~\cite{lee2018primal}. For the comparative method, the discount factor, the approximation factor, and the initial value number are set as $\alpha=1$, $M=1000$, and $r=4$, respectively. Also, the initial value vectors are designed as 
	\begin{equation}\nonumber
		\  V=\begin{bmatrix} v_1^T \\ v_2^T \\v_3^T\\v_4^T \end{bmatrix} =	\begin{bmatrix} 1 & 1 & 1\\1 & 10 & 2\\3 & 2 & 0\\5 & 2 & 0 \end{bmatrix}.
	\end{equation}
	The definitions of the above-mentioned parameters and further details of the comparative method can be found in~\cite{lee2018primal}. Considering the structural constraints in this problem, the projection is applied in the comparative method. Besides, as required by the primal-dual Q-learning algorithm, the continuous system is transformed into its corresponding discrete system with a sampling time of 0.01 s. For the sake of fairness, the iteration number is also set as 10, the initial values of the controller gain and the state variables are set as the same as in our proposed method, i.e. $K_0 =\begin{bmatrix} 0&-5 \end{bmatrix}$ and $x_0 =\begin{bmatrix} 1&0 \end{bmatrix}^T$. After 10 iterations, the cost is given by $1.3251$ and the controller gain is given by \begin{equation} \nonumber K^*=\begin{bmatrix}  0  &  -3.2311\end{bmatrix}.\end{equation} With this  gain matrix, the system response is shown as the dashed lines in Fig.~\ref{fig:response1}. As shown in the figure, it is illustrated that the performance attained by our proposed method and the primal-dual Q-learning algorithm is rather close. According to the comparison of the cost, our approach leads to a slightly lower one than that of the primal-dual Q-learning algorithm, and this illustrates the correctness of our approach.}
	

	\begin{example}\label{exam:2}
Consider the nominal model of a linear 4th-order system, where
		\begin{align}\nonumber
			\dot {{x}}&=A {x} +B u, \quad x(0)=x_0,\\
			u&=Kx,\nonumber
		\end{align}
		and
		\begin{equation}\nonumber
		A =
		\begin{bmatrix}
		0 &  1 & 0 & 0\\
		12& 0 & -10 & 0 \\
		0 & 0 & 0 & 1\\
		-12 & 0 & 25 & 0
		\end{bmatrix}, \quad 	  B  =
		\begin{bmatrix}
		0 & 0 \\
		1 & -2\\
		0 & 0\\
		-2 & 5
		\end{bmatrix}.\end{equation}
\textcolor{black}{Assume the system is perturbed by model uncertainties (under which is considered to be operated in the actual condition), and all the non-zero elements in the state matrix $A$ varies within $\pm20\%$ of their nominal values, respectively. In this case, we also consider such a perturbed system, where
\begin{align}\nonumber
			\dot {{x}}&=\tilde A {x} +B u, \quad x(0)=x_0,\\
			u&=Kx,\nonumber
		\end{align}
		and
		\begin{equation}\nonumber
		\tilde A =
		\begin{bmatrix}
		0 &  1 & 0 & 0\\
		9.8& 0 & -9.8 & 0 \\
		0 & 0 & 0 & 1\\
		-9.8 & 0 & 29.4 & 0
		\end{bmatrix}, \quad 	  B  =
		\begin{bmatrix}
		0 & 0 \\
		1 & -2\\
		0 & 0\\
		-2 & 5
		\end{bmatrix}.\end{equation}
		}
		
In terms of the actual operating condition where the system is perturbed, we aim to determine an optimal controller \begin{equation}\nonumber K=\begin{bmatrix} k_1 & k_2 & k_3 & k_4 \\k_5 & k_6 & k_7 & k_8 \end{bmatrix},\end{equation} with $k_3=k_4=k_5=k_6=0$ (which leads to a decentralized control system), such that
		\begin{equation}\nonumber
		\mathop{\textup{min}}\limits_{K} \hspace{1.5mm} \int_0^\infty \left(x^TQx+u^TRu \right) \,dt,
		\end{equation}
		where  $Q=\begin{bmatrix} 1 & 0 & 0 & 0\\0 & 0 & 0 & 0 \\0 & 0 & 1 &0\\ 0 & 0 & 0 & 0\end{bmatrix}$, $R=\begin{bmatrix}  1 & 0 \\0 & 1\end{bmatrix}$.
	\end{example}
The structural constraints on $K$ can be stated as
	\begin{gather}
		G_1KH_1 =0, \quad G_1=\begin{bmatrix} 1 & 0   \end{bmatrix}, \quad H_1=\begin{bmatrix}   0 & 0 & 1 & 0 \end{bmatrix}^T; \nonumber\\
		G_2KH_2 =0, \quad G_2=\begin{bmatrix} 1 & 0   \end{bmatrix}, \quad H_2=\begin{bmatrix}   0 & 0 & 0 & 1 \end{bmatrix}^T; \nonumber\\
		G_3KH_3 =0, \quad G_3=\begin{bmatrix} 0 & 1   \end{bmatrix}, \quad H_3=\begin{bmatrix}   1 & 0 & 0 & 0 \end{bmatrix}^T; \nonumber\\
		G_4KH_4 =0, \quad G_4=\begin{bmatrix} 0 & 1   \end{bmatrix}, \quad H_4=\begin{bmatrix}   0 & 1 & 0 & 0 \end{bmatrix}^T. \nonumber\nonumber
	\end{gather} 

\begin{figure}
		\centering
		\includegraphics[trim=0 0 0 0,width=0.9\columnwidth]{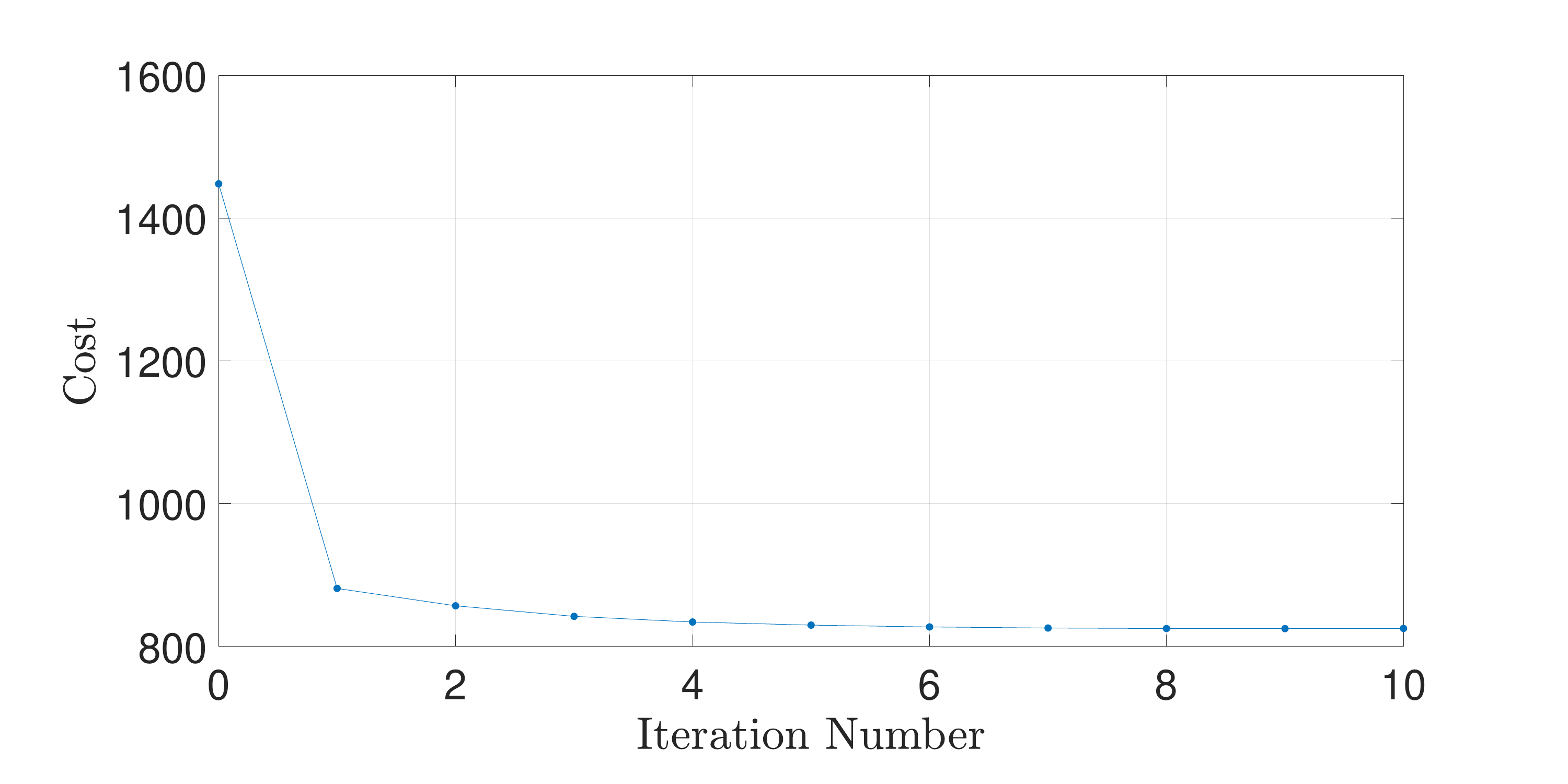}
		\caption{Change of the cost during iterations in Example 2.}
		\label{fig:cost2}
	\end{figure}

\begin{figure}
		\centering
		\includegraphics[trim=0 0 0 0,width=0.9\columnwidth]{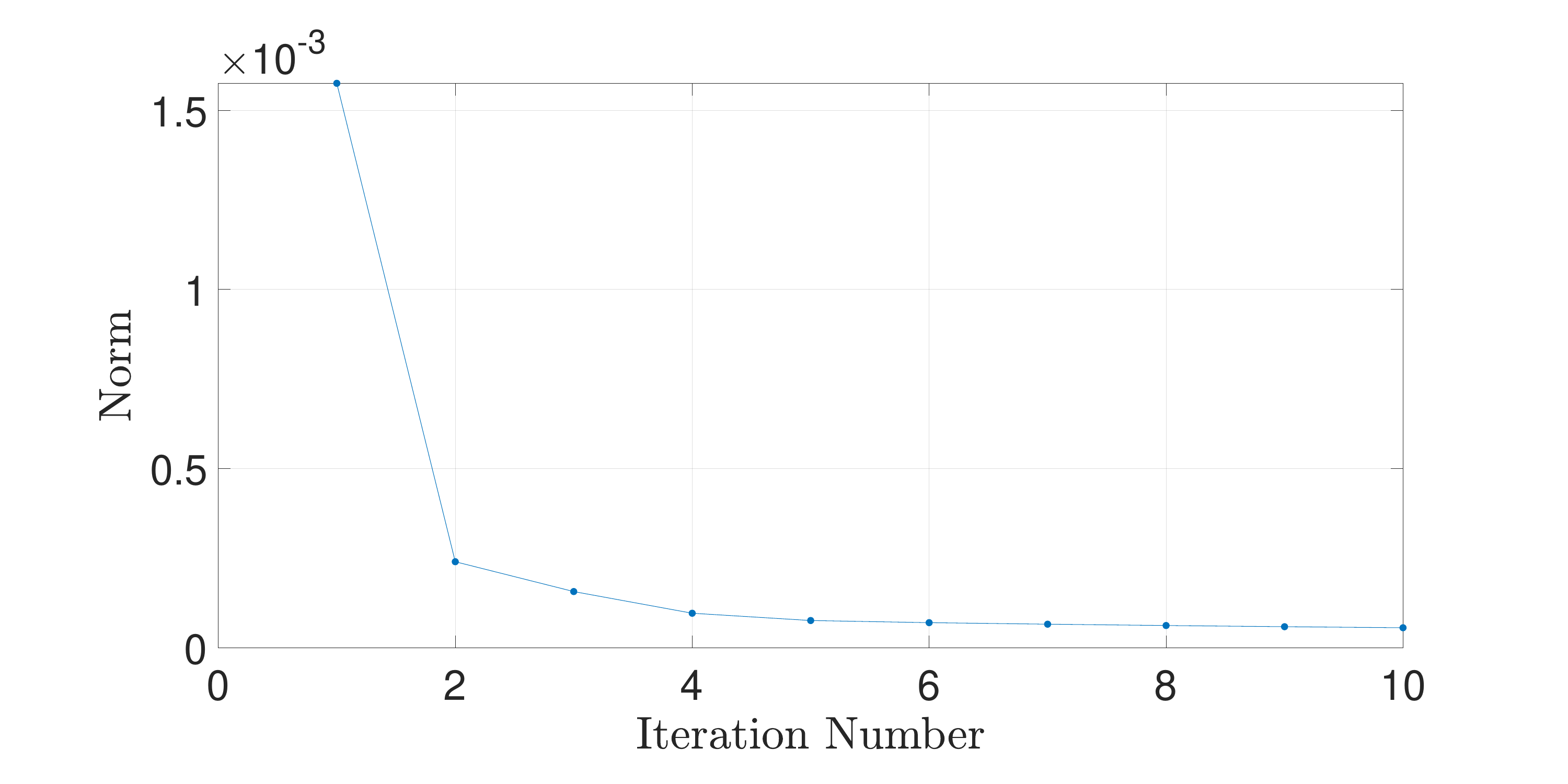}
		\caption{Change of the norm during iterations in Example 2.}
		\label{fig:norm2}
	\end{figure}

\begin{figure}[t]
		\centering
		\includegraphics[trim=0 0 0 0,width=0.9\columnwidth]{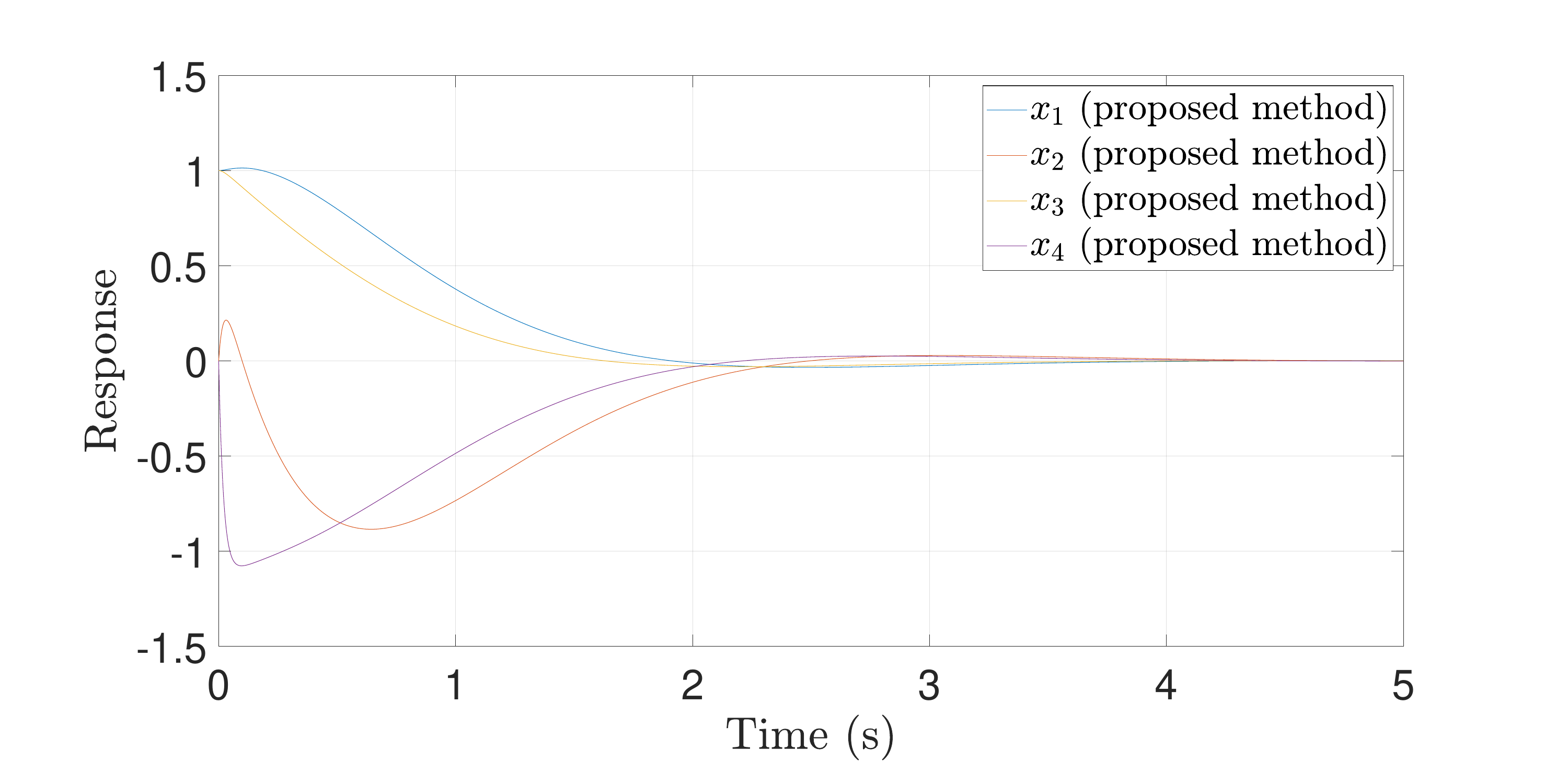}
		\caption{System response of the proposed method in Example 2.}
		\label{fig:response2}
	\end{figure}

In this example, the initialized controller gain is set as \begin{equation}\nonumber K=\begin{bmatrix} -50&-20 & 0 &0\\0&0& -20&-6 \end{bmatrix},\end{equation} which stabilizes the closed-loop system with a cost of $1.4480 \times 10^{3}$. Besides, the iteration number is set as 10 and the initial condition is $x_0=\begin{bmatrix}
	1&0&1&0
	\end{bmatrix}^T$. In the proposed data-driven method, Fig.~\ref{fig:cost2} and Fig.~\ref{fig:norm2} show the changes of the cost and the norm during 10 iterations in the optimization. After 10 iterations, the cost is given by $8.2534\times 10^{2}$ with a norm of $5.6196 \times 10^{-5}$, and the optimal controller gain is given by \begin{equation}\nonumber K^*=\begin{bmatrix} -57.8053  & -18.8322   & 0 &0\\0&0&  -38.0592 &  -7.3845 \end{bmatrix}.\end{equation}
	With $x_0=\begin{bmatrix}
	1&0&1&0
	\end{bmatrix}^T$ and the optimal gain matrix, the system response is shown in Fig.~\ref{fig:response2}.

\begin{figure}[t]
		\centering
		\includegraphics[trim=0 0 0 0,width=0.9\columnwidth]{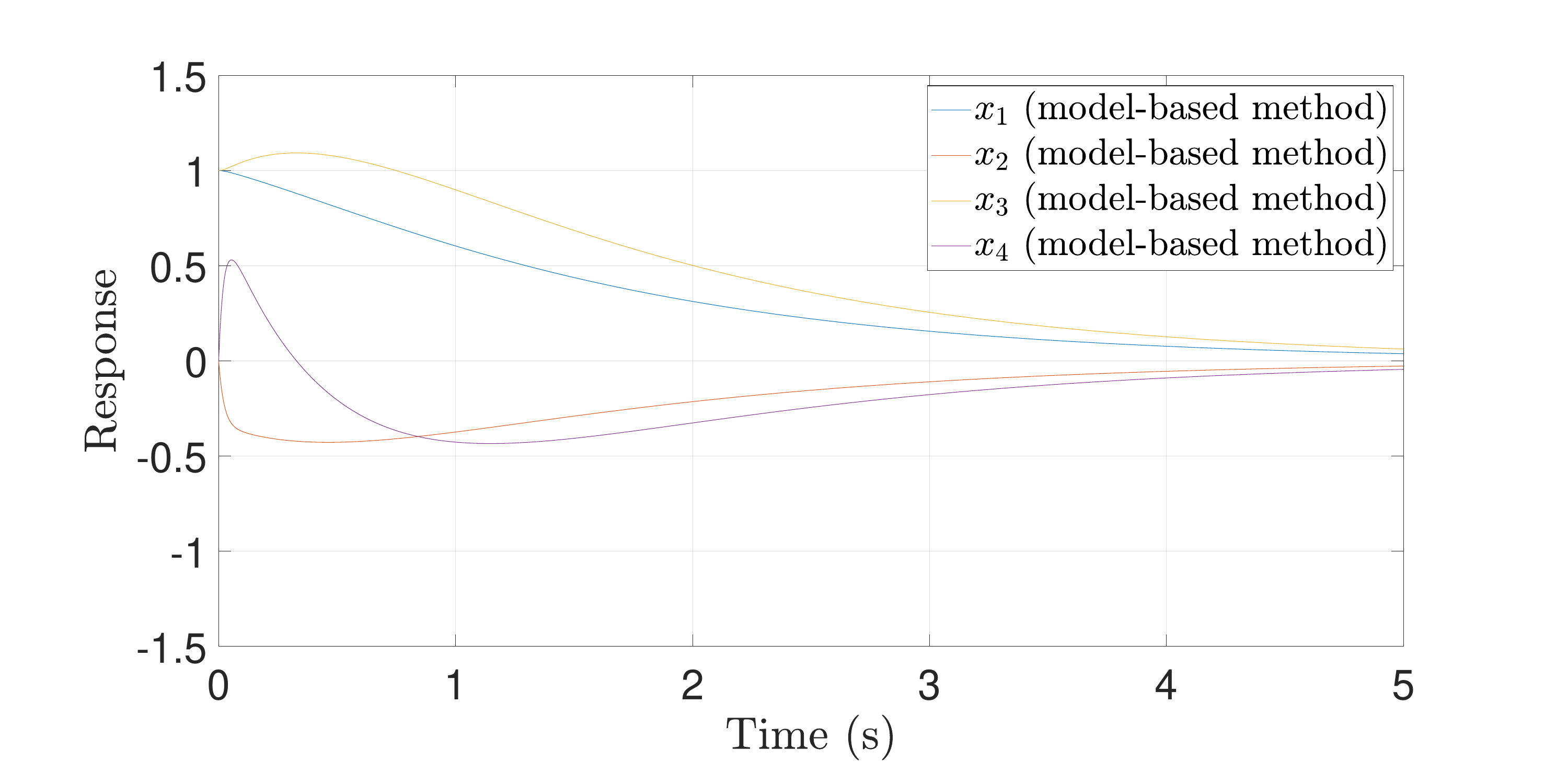}
		\caption{\textcolor{black}{System response of the comparative method in Example 2.}}
		\label{fig:response2compare}
	\end{figure}

\textcolor{black}{To further demonstrate the effectiveness and superiority of the proposed method, a model-based constrained linear quadratic optimization method~\cite{geromel1982optimal} is used in this study as the comparative method. It is worthwhile to note that, for the model-based constrained linear quadratic optimization method, the knowledge of the system model is needed when calculating the controller gain. Thus, in this case study, the nominal system is used to calculate the result of the controller gain, which is further applied to the actual perturbed system to observe the performance. With the comparative model-based method, all the settings regarding the initial gain and iteration number remain the same as those in the proposed data-driven method for the sake of fairness. After 10 iterations, the result of the controller gain is given by \begin{equation}\nonumber K^*=\begin{bmatrix} -54.0897  & -17.8872   & 0 &0\\0&0&  -19.2366 &  -6.4998 \end{bmatrix},\end{equation}  and the cost is given by $1.2060\times 10^{3}$ with a norm of $3.0566 \times 10^{-3}$. With the determined controller gain, the system response is illustrated in Fig.~\ref{fig:response2compare}. }


	\textcolor{black}{
	    Comparing the two methods, it can be easily seen that the cost with the proposed data-driven method is significantly lower than the one with the comparative model-based method, and the performance of the proposed data-driven method is obviously better than the one of the comparative model-based method. The reasonable explanation of the phenomena is that, the optimum of the result highly relies on the system model with the model-based linear quadratic optimization methods; while with the data-driven approach, the optimal solution is attained only from the data rather than the system model. As a result, the superiority over the model-based approach is demonstrated over the model-based approach in situations of significant model uncertainties.}


	\section{Conclusion}
	In this work, a data-driven optimization algorithm is developed to solve the LQR problem with structural constraints.
	The gradient of the objective function with respect to the gain matrix is derived for both the single-input and the multi-input cases,
	and then it is projected onto the related hyperplanes characterizing the zero-element constraints.
	In this way, the gain matrix is iteratively updated towards the optimum with structural constraints preserved.
	The proposed algorithm provides a suitably general and effective methodology
	for solving the constrained linear quadratic optimization problem
	without the need for precise modeling.
	In the work here, the formulation has also been evaluated and verified on numerical examples, where the effectiveness of the theoretical results are successfully and clearly demonstrated. Possible future works include the generalization of the results to nonlinear systems, as well as the problems with cost functions of other types (such as the non-quadratic cost function). 	And certainly for the
		pertinent observations that relate to the important aspect on further improvement of the practical scalability of the proposed method,
		it can be stated here that this is a significant and substantial aspect
		where for this approach of the proposed methodology,
		the computations involved in
		the control of systems with larger dimensions 
		will also be an important matter
		as part of the future work.

\bibliographystyle{IEEEtran}
\bibliography{IEEEabrv,Reference}

\end{document}